\documentclass[10pt,final,journal]{IEEEtran}
\usepackage{amsmath}
\usepackage{bm}
\usepackage{amssymb}
\usepackage{amsthm}
\usepackage{graphicx}
\usepackage{bbold}
\newcommand{\Kf}{K\hspace{-0.8mm}f}
\usepackage[colorlinks, bookmarks=false, citecolor=blue, linkcolor=red, urlcolor=blue]{hyperref}

\newtheorem{theorem}{Theorem}

\newtheorem{proposition}{Proposition}
\newtheorem{assumption}{Assumption}
\newtheorem{remark}{Remark}
\newtheorem{conjecture}{Conjecture}

\begin{document}
\title{Optimal placement of inertia and primary control : a matrix perturbation theory approach}
\author{Laurent~Pagnier,~\IEEEmembership{Member,~IEEE} and Philippe~Jacquod,~\IEEEmembership{Member,~IEEE}}
\maketitle
\begin{abstract}
The increasing penetration of inertialess new renewable energy sources 
reduces the overall mechanical inertia available in power grids and accordingly
raises a number of issues of  grid stability over short to medium time scales. 
It has been suggested that
this reduction of overall inertia can be compensated to some extent by the deployment of 
substitution inertia - synthetic inertia, flywheels or synchronous condensers.
Of particular importance is to optimize the placement of the limited available substitution inertia,
to mitigate voltage angle and frequency disturbances 
following a fault such as an abrupt power loss. Performance measures in the form of ${\cal H}_2-$norms have been recently
introduced to evaluate the overall magnitude of such disturbances on an electric power 
grid. However, despite the mathematical conveniance
of these measures, analytical results can be obtained only under rather restrictive assumptions of 
uniform damping ratio, or homogeneous distribution of inertia and/or primary control in the 
system. Here, we introduce matrix perturbation theory  to obtain analytical results for optimal inertia and 
primary control placement where both are heterogeneous. 
Armed with that efficient tool, we construct two simple algorithms that independently determine the optimal geographical 
distribution of inertia and primary control. These algorithms are then implemented on a model 
of the synchronous transmission grid of continental Europe. 
We find that the optimal distribution of inertia is geographically homogeneous but that primary
control should be mainly located on the slow modes of the network, where the intrinsic grid
dynamics takes more time to damp frequency disturbances.  
\end{abstract}

\section{Introduction}

The penetration of new renewable energy sources (RES) such as photovoltaic panels and wind turbines is increasing
in most electric power grids worldwide. In their current configuration, these energy sources are essentially inertialess,
and their increased penetration leads to low inertia situations in periods of high RES production~\cite{ulbig2014impact}. 
This raises important 
issues of power grid stability, which is of higher concern to transmission system operators than the volatility of 
the RES productions~\cite{Mil15,Win15}. The substitution of traditional productions based on
synchronous machines with inertialess RES may in particular lead to 
geographically inhomogeneous inertia profiles.  
It has been suggested to deploy substitution inertia -- synthetic inertia, flywheels or synchronous condensers -- 
to compensate locally or globally for missing inertia. Two related question naturally arise, which are (i) where is it safe to 
substitute synchronous machines with inertialess RES and (ii) where is it optimal to distribute substitution inertia ?
Problem (ii) has been investigated in small power grid models with up to a dozen 
buses, optimizing the geographical distribution of inertia against cost functions based 
on eigenvalue damping ratios~\cite{borsche2015effects}, $\mathcal{H}_p$-norms~\cite{mevsanovic2016comparison,poolla2017optimal} or RoCoF~\cite{paganini2017global,bor18} and frequency
excursions~\cite{bor18}.
Investigations of problem (i) on large power grids 
emphasized the importance of the geographical extent of the slow network modes~\cite{pagnier2019inertia}.
Numerical optimization can certainly be performed for any given network on a case-by-case basis, however
it is highly desirable to shed light on the problem with analytical results. So far, such results have been either 
restricted to small systems or derived under the assumption that damping and inertia parameters,
or their ratio, are homogeneous. In this manuscript we go beyond these assumptions and construct an
approach applicable to large power grids with inhomogeneous independent damping and inertia parameters. 

Inspired by theoretical physics, we introduce {\it matrix perturbation theory}~\cite{stewart1990matrix} as an
analytical tool to tackle this  problem. That method is widely used in quantum 
physics, where it delivers approximate solutions
to complex, perturbed 
problems, extrapolated from known, exact solutions of integrable problems~\cite{macintyre}. 
The approximation is valid as long as the difference between the two problems is small and it makes sense to 
consider the full, complex problem as a {\it perturbation} of the exactly soluble, simpler problem. 
The procedure is spectral in nature. It identifies a small, dimensionless parameter in which 
eigenvalues and eigenvectors of the perturbed problem can be systematically expanded in a series not unlike a Taylor-expansion
about the unperturbed, integrable solution. 
Depending on the value of that small parameter, the series can be truncated at low orders already and still deliver 
rather accurate results. In the context of electric power grids,
the method was applied for instance in Ref.~\cite{coletta2018performance}, where quadratic performance measures similar to those discussed below were calculated following a line fault, starting from 
the eigenvalues and eigenvectors of the network Laplacian before the fault. 

In this paper, we apply matrix perturbation theory~\cite{stewart1990matrix} to calculate performance measures
following an abrupt power loss in the case of a transmission power grid with geographically  
inhomogeneous inertia and damping (primary control) parameters. Our perturbation theory is an expansion in two parameters
which are 
the maximal deviations $\delta m$ and $\delta d$ of the rotational inertia and damping parameters from their
average values $m$ and $d$. The approach is valid as long as these local deviations are small, 
$|\,\delta m/m\,|<1$, $|\,\delta d/d\,|<1$. These conditions tolerate in principle that inertia and damping parameters
vanish or are twice as large as their average on some buses. 
The main step forward brought about by our approach is that we are able 
to derive analytical results without relying on the often used homogeneity assumptions that damping, inertia or
their ratio is constant - assumptions which are not satisfied in real electric power grids. 
Our main results are given in Theorems~\ref{thm:lin_opt_inertia} and \ref{thm:lin_opt_control} below, which formulate
algorithms for optimal placement of local inertia and damping parameters.
The spectral decomposition approach used here has recently drawn the attention of a number of groups and has been used
to calculate performance measures in power grids and consensus algorithms e.g. in \cite{paganini2017global,Pir17,Guo18,Por18}. 

The article is organized as follows. Section~\ref{section:homogeneous} deals with the case where inertia and primary control are 
uniformly distributed in the system. The performance measure that quantifies system disturbances is introduced 
and we calculate its value for abrupt power losses. In Section~\ref{section:perturbation} we apply matrix perturbation theory to  calculate the sensitivities of our measure in local variations of inertia and primary control.  Section~\ref{section:optimal} presents the optional placement of inertia and primary control in the case of weak inhomogeneity. In Section~\ref{section:applications} we apply our optimal placements to the continental European grid.  Section~\ref{section:conclusion} concludes our article.

\section{Homogeneous case}\label{section:homogeneous}

We are interested in the dynamical response of an electric power grid to a disturbance such as an abrupt power loss. To that end,
we consider the power system dynamics in the lossless line approximation, which is a standard approximation used for 
high voltage transmission grids~\cite{machowski2008power}. That dynamics is governed by the swing equations,
\begin{equation}\label{eq:swing}
m_i\dot\omega_i+d_i\omega_i=P_i-\sum_jB_{ij}\sin(\theta_i-\theta_j)\, ,
\end{equation}
which determine the time-evolution of the voltage angles $\theta_i$ and frequencies $\omega_i=\dot{\theta}_i$
at each of the $N$ buses labelled 
$i$ in the power grid in a rotating frame such that $\omega_i$ measures the angle frequency deviation 
to the rated grid frequency
of 50 or 60 Hz. Each bus is characterized by inertia, $m_i$, and damping, $d_i$, parameters
and $P_i$ is the active power injected ($P_i>0$) or extracted ($P_i<0$) at bus $i$. We introduce the 
damping ratio $\gamma_i \equiv d_i/m_i$.
Buses are connected
to one another via lines with susceptances $B_{ij}$. Stationary solutions $\{\theta_i^{(0)}\}$ are power flow solutions 
determined by $P_i=\sum_jB_{ij}\sin(\theta_i^{(0)}-\theta_j^{(0)})$. 
Under a change in active power $P_i \rightarrow P_i + \delta P_i$, linearizing the dynamics 
about such a solution with $\theta_i(t) = \theta_i^{(0)}+\delta \theta_i(t)$ gives, in matrix form, 
\begin{equation}
\bm M\bm{\dot\omega}+\bm D \bm \omega= \bm{\delta P} - \bm{L \delta \theta} \,,\label{eq:swinglin}
\end{equation}
where $\bm M={\rm diag}(\{m_i\})$, $\bm D={\rm diag}(\{d_i\})$ and voltage angles and frequencies are cast into vectors
$\bm{\delta\theta}$ and $\bm \omega\equiv\bm{\delta\dot\theta}$. 
The Laplacian matrix $\bm L$ has matrix elements
$L_{ij} = -B_{ij} \cos(\theta_i^{(0)}-\theta_j^{(0)})$, for $i \ne j$ and 
$L_{ii} = \sum_{k} B_{ik} \cos(\theta_i^{(0)}-\theta_k^{(0)})$. 

\subsection{Exact solution for homogeneous damping ratio}

When the damping ratio is constant, $d_i/m_i=\gamma_i = \gamma$, $\forall i$, \eqref{eq:swinglin} can be integrated
exactly~\cite{paganini2017global,coletta2018performance}. To see this we first 
transform angle coordinates as
$\bm{\delta\theta}=\bm{M}^{-1/2}\bm{\delta\theta_M}$ to obtain
\begin{equation}
\bm{\dot{\omega}_M}+\underbrace{\bm M^{-1} \bm D}_{\bm \Gamma}\bm{\omega_M}+\underbrace{\bm M^{-1/2}\bm L \bm M^{-1/2}}_{\bm{L_M}}\bm{\delta\theta_M}=\bm M^{-1/2}\bm{\delta P}\,,\label{eq:omegaM}
\end{equation}
where we introduced the diagonal matrix $\bm \Gamma = {\rm diag}(\{d_i/m_i\}) \equiv {\rm diag}(\{\gamma_i\})$. 
The inertia-weighted matrix $\bm{L_M}$ is real and symmetric, therefore it can be diagonalized 
\begin{equation}\label{eq:lm}
\bm{L_M} = \bm U^{\top}\bm\Lambda \bm U
\end{equation}
with an orthogonal matrix $\bm U$, the $\alpha^{\rm th}$ row of which gives the components $u_{\alpha,i}$, $i=1, \ldots N$ 
of the $\alpha^{\rm th}$ eigenvector $\bm u_\alpha$ of $\bm{L_M}$. The diagonal matrix 
$\bm \Lambda={\rm diag}(\{\lambda_1=0,\lambda_2,\cdots,\lambda_N\})$ contains the eigenvalues of $\bm L_M$
with $\lambda_\alpha<\lambda_{\alpha+1}$. For connected networks only the smallest eigenvalue $\lambda_1$ vanishes, 
which follows from the zero row and column sum property of the Laplacian
matrix $\bm{L_M}$. Rewriting \eqref{eq:omegaM} in the basis diagonalizing $\bm L_M$ gives
\begin{equation}
\bm{\ddot\xi}+\bm U \bm\Gamma \bm U^\top \bm{\dot\xi} +\bm\Lambda \bm\xi=\bm U\bm M^{-1/2}\bm{\delta P}\, ,\label{eq:xi_ode}
\end{equation}
where $\bm{\delta\theta_M}=\bm U^\top \bm \xi$. This change of coordinates is nothing but a spectral decomposition of angle
deviations $\bm{\delta\theta_M}$ into their components in the basis of eigenvectors of $\bm L_{\bm M}$. These components
are cast in the vector $\bm \xi$. 
The formulation \eqref{eq:xi_ode} of the problem makes it clearer that, if $\bm \Gamma$ 
is a multiple of identity, the problem can be recast as a diagonal ordinary differential equation problem 
that can be exactly integrated. This is done below in \eqref{eq:chi_ode}, and 
provides an exact solution about which we will construct a matrix perturbation theory in the next sections.  

\begin{proposition}[Unperturbed evolution]  \label{proposition1}
For an abrupt power loss,
$\bm{\delta P}(t) = \bm{\delta P} \, \Theta(t)$, with the Heaviside step function defined by $\Theta(t>0)=1$, $\Theta(t<0)=0$,
and with homogeneous damping ratio, 
$\bm \Gamma=\gamma \, \mathbb{1}$ with the $N \times N$ identity matrix $\mathbb{1}$, the frequency 
coordinates $\dot\xi_\alpha$ evolve independently as
\begin{align}
\dot\xi_\alpha(t)&=\frac{2\mathcal{P}_\alpha}{f_\alpha }e^{-\gamma t/2}\sin\Big(\frac{f_\alpha t}{2}\Big)\, \text{, $\forall \alpha >1$,}\label{eq:dchii}
\end{align}
where $f_\alpha =\sqrt{4\lambda_\alpha -\gamma^2}$ and $\mathcal{P}_\alpha=\sum_i u_{\alpha i}\,\delta P_i / m_i^{1/2}$.
\end{proposition} 
This result generalizes Theorem III.3 of \cite{Guo18}.
\begin{IEEEproof}
The proof goes along the lines of the diagonalization procedure proposed in 
\cite{paganini2017global,coletta2018performance,coletta2018transienta}. Equation
\eqref{eq:xi_ode} can be rewritten as
\begin{equation}
\frac{{\rm d}}{{\rm d}t}\left[\!\begin{array}{c}
\bm\xi\\
\bm{\dot\xi}
\end{array}\!\right]=\underbrace{
\left[\!\begin{array}{cc}
\mathbb{0}_{N\times N}&\!\!\!\!\mathbb{1}\\
\!\!\!\!-\bm\Lambda&\!\!\!\!\!\!-\gamma \, \mathbb{1}
\end{array}\!\right]}_{\bm H_0}\!\left[\!\begin{array}{c}
\bm\xi\\
\bm{\dot\xi}
\end{array}\!\right]+
\bigg[\!\begin{array}{c}
\mathbb{0}_{N\times 1} \\
\bm{\mathcal{P}}
\end{array}\!\bigg]\,,\label{eq:H}
\end{equation}
where $\bm{\mathcal{P}}=\bm U\bm M^{-1/2} \bm{\delta P}$ and $\mathbb{0}_{N\times M}$ is the $N\times M$ matrix of zeroes. The matrix 
$\bm H_0$ is block-diagonal up to a permutation of rows and columns \cite{coletta2018performance}, and can easily be diagonalized 
block by block, where each $2 \times 2$ block corresponds to one of the eigenvalues $\lambda_\alpha$ of $\bm{L_M}$. 
The $\alpha^{\rm th}$ block is diagonalized by the transformation
\begin{align}
\left[\!\begin{array}{c}
\chi^{(0)}_{\alpha+}\\
\chi^{(0)}_{\alpha-}
\end{array}\!\right]
&=\bm{T}_{\alpha}^{L}\left[\!\begin{array}{c}
\xi_\alpha\\
\dot\xi_\alpha
\end{array}\!\right]\,, \;\;\; \bm{T}_{\alpha}^{L}\equiv 
\frac{i}{f_\alpha }\left[\!\begin{array}{cc}
\phantom{-}\mu_{\alpha-}^{(0)} & -1\\
-\mu_{\alpha+}^{(0)} &\phantom{-}1
\end{array}\!\right] \, ,\label{eq:TiL}
\\
\left[\!\begin{array}{c}
\xi_\alpha\\
\dot\xi_\alpha
\end{array}\!\right]
&= \bm{T}_{\alpha}^{R}\left[\!\begin{array}{c}
\chi^{(0)}_{\alpha +}\\
\chi^{(0)}_{\alpha -}
\end{array}\!\right]\, , \;\;\;
\bm{T}_{\alpha}^{R} \equiv
\left[\!\begin{array}{cc}
1 & 1\\
\mu_{\alpha+}^{(0)} & \mu_{\alpha-}^{(0)} 
\end{array}\!\right]\, ,\label{eq:TiR}
\end{align}
with the eigenvalues $\mu_{\alpha\pm}^{(0)}$ of the $\alpha^{\rm th}$ block,
\begin{equation}
\mu_{\alpha\pm}^{(0)}=-\frac{1}{2}(\gamma \mp if_\alpha)\, .
\end{equation}
The two rows (columns) of $\bm{T}_{\alpha}^{L}$ ($\bm{T}_{\alpha}^{R}$) give the nonzero components  of the 
two left (right) eigenvectors $\bm t_{\alpha \pm} ^{(0)L}$ ($\bm t_{\alpha \pm}^{(0)R}$) of $\bm H_0$.
Following this transformation, \eqref{eq:H} reads
\begin{equation}
\frac{{\rm d}}{{\rm d}t}\left[\!\!\begin{array}{c}
\chi^{(0)}_{\alpha +}\\
\chi^{(0)}_{\alpha -}
\end{array}\!\!\right]=
\left[\!\!\begin{array}{cc}
\mu_{\alpha +}^{(0)} & 0\\
0 & \mu_{\alpha -}^{(0)}
\end{array}\!\!\right]\!\left[\!\!\begin{array}{c}
\chi^{(0)}_{\alpha +}\\
\chi^{(0)}_{\alpha -}
\end{array}\!\!\right]+
\frac{i}{f_{\alpha}}\left[\!\!\begin{array}{c}
-\mathcal{P}_\alpha \\
\phantom{-}\mathcal{P}_\alpha 
\end{array}\!\!\right]\,.\label{eq:chi_ode}
\end{equation}
The solutions of \eqref{eq:chi_ode} are
\begin{align}
\chi^{(0)}_{\alpha\pm}&=\pm\frac{i\,\mathcal{P}_\alpha }{f_\alpha \mu_{\alpha\pm}^{(0)}}\Big(1-e^{\mu_{\alpha\pm}^{(0)} t}\Big)\,,\;\forall \alpha>1\label{eq:chiipm} \, .
\end{align}
Inserting \eqref{eq:chiipm} back into \eqref{eq:TiL}, one finally finds \eqref{eq:dchii}
which proves the proposition.
\end{IEEEproof}

\subsection{Performance measure}
We want to mitigate disturbances following an abrupt power loss. To that end, we use performance measures
which evaluate the overall disturbance magnitude over time and the whole power grid.
Performance measures have been proposed, which
can be formulated as $\mathcal{L}_2$ and squared $\mathcal{H}_2$ norms of 
linear   
systems~\cite{poolla2017optimal,paganini2017global,coletta2018performance,coletta2018transienta,tegling2015price,Fardad14,gayme16,siami2016fundamental,tyloo2018robustness,guo2019performance} and are time-integrated quadratic forms in the angle, 
$\bm{\delta \theta}$, or frequency, $\bm \omega$, deviations. Here we focus on frequency deviations and use the following 
performance measure
\begin{equation}
\mathcal{M}=\intop_0^\infty\big(\bm \omega^\top -\bm{\bar\omega}^\top \big)\bm M\big(\bm \omega-\bm{\bar\omega}\big){\rm d}t\,,\label{eq:measure}
\end{equation}
where $\bm{\bar\omega}= (\omega_{\rm sys}, \omega_{\rm sys}, \ldots\omega_{\rm sys})^\top$ is the 
instantaneous average frequency vector with components
\begin{equation}
\omega_{\rm sys}(t) = \sum_i m_i \omega_i(t) \Big / \sum_i m_i \, .
\end{equation}
It is straightforward to see that $\mathcal{M}$ reads
\begin{equation}\label{eq:M}
\mathcal{M}=\intop_0^\infty\sum_{\alpha>1}\dot\xi_\alpha^2(t){\rm d}t\, ,
\end{equation}
when rewritten in the eigenbasis of $\bm L_{\bm M}$, once one notices that the first eigenvector of $\bm L_{\bm M}$
(the one with zero eigenvalue) has components $u_{1i}=\sqrt{m_i}/\sqrt{\sum_j m_j}$.
\begin{proposition}
For an abrupt power loss,
$\bm{\delta P}(t) = \bm{\delta P} \, \Theta(t)$ on a single bus labeled $b$, 
$\delta P_i=\delta_{ib} \, \delta P$,
and with an homogeneous damping ratio, 
$\bm \Gamma=\gamma \, \mathbb{1}$ with the $N \times N$ identity matrix $\mathbb{1}$, 
\begin{equation}
\mathcal{M}_b=\frac{\delta P^2}{2\gamma m_b}\sum_{\alpha>1}\frac{u_{\alpha b}^2}{\lambda_\alpha }\,,\label{eq:Mb}
\end{equation}
in terms of the eigenvalues $\lambda_\alpha$ and the components $u_{\alpha b}$ of the eigenvectors
$\bm u_\alpha$ of $\bm L_{\bm M}$.
\end{proposition}
Note that we introduced the subscript $b$ to indicate that the fault is localized on that bus only. 
The power loss is modeled as $P_i = P_i^{(0)} -\delta P_i \, \Theta(t)$ with  $\delta P_i = \delta_{ib} \, \delta P $ with the Kronecker
symbol $\delta_{ib}=1$ if $i = b$ and 0 otherwise. 
\begin{IEEEproof}
Equation  \eqref{eq:dchii} straightforwardly gives
\begin{equation}
\intop_0^\infty\dot\xi_\alpha^2(t){\rm d}t = \frac{\delta P^2\,u_{\alpha b}^2}{2\gamma \, m_b \, \lambda_\alpha }\,, \alpha>1\, ,
\end{equation}
which, when summed over $\alpha >1$ gives \eqref{eq:Mb}.
\end{IEEEproof}
\begin{remark}
For homogeneous inertia coefficients, $\bm M=m\mathbb{1}$, the eigenvectors and eigenvalues 
of the inertia-weighted Laplacian $\bm L_{\bm M}$
defined in \eqref{eq:omegaM} are given by
$\bm u_{\alpha}=\bm u_{\alpha}^{(0)}$, and 
$\lambda_\alpha=m^{-1}\lambda_\alpha^{(0)}$, in terms of the eigenvectors $\bm u_{\alpha}^{(0)}$ and eigenvalues
$\lambda_\alpha^{(0)}$ of the Laplacian $\bm L$. In that case, the performance measure reads
\begin{equation}
\mathcal{M}_b^{(0)}=\frac{\delta P^2}{2\gamma}\sum_{\alpha>1}\frac{u_{\alpha b}^{(0)2}}{\lambda_\alpha ^{(0)}}\,,
\label{eq:crit}
\end{equation}
where the superscript $^{(0)}$ refers to inertia homogeneity. This expression has an interesting graph theoretic 
interpreation. We recall the definitions of the resistance distance $\Omega_{ij}$ between two nodes on the network,
the associated centrality $C_j$ and the generalized Kirchhoff indices $\Kf_p$~\cite{tyloo2018robustness,Kle93},
\begin{eqnarray}
\Omega_{ij} &=& \bm L_{ii}^\dagger + \bm L_{jj}^\dagger -\bm L_{ij}^\dagger -\bm L_{ji}^\dagger \, , \label{eq:rdistance} \\
C_j &=& N \, \Big(\sum_{i}\Omega_{ij}\Big)^{-1} \, , \label{eq:centrality} \\
\Kf_p &=& N \sum_{\alpha>1} \lambda_\alpha^{-p} \, ,
\label{eq:kfp} 
\end{eqnarray}
where $\bm L^\dagger$ is the Moore--Penrose pseudo inverse of $\bm L$. 
With these definitions, one can show that \cite{Por18,tyloo2018robustness,tyloo2018key}
\begin{equation}\label{eq:kf}
\sum_{\alpha>1} \frac{u_{\alpha b}^{(0)2}}{\lambda_\alpha ^{(0)}} = C_b^{-1} - N^{-2} \Kf_1 \, ,
\end{equation}
by using the spectral representation of the resistance distance~\cite{Gut96,Zhu96}
\begin{equation}
\Omega_{ib}=\sum_{\alpha>1}\big(u_{\alpha i}^{(0)}-u_{\alpha b}^{(0)}\big)^2\big/\lambda_{\alpha}^{(0)}\,.\label{eq:res_dist}
\end{equation}
Because $\Kf_1$ is a global quantity
characterizing the network, it follows from \eqref{eq:crit} with \eqref{eq:kf} that, 
when inertia and primary control are homogeneously distributed in the system, the disturbance magnitude as measured
by $\mathcal{M}_b^{(0)}$ is larger for disturbances on peripheral nodes~\cite{pagnier2019inertia,Gam17}. 
\end{remark}

\section{Matrix perturbation}\label{section:perturbation}

The previous section treats the case where inertia and primary control are uniformly distributed in the system. Our goal is to lift that restriction and to obtain $\mathcal{M}_b$ when some mild inhomogeneities are present. We parametrize these inhomogeneities by writing
\begin{align}
m_i&= m+\delta m\, r_i\,,\label{eq:mi}\\
d_i&= m_i\gamma_i=(m+\delta m\, r_i)(\gamma +\delta\gamma\, a_i)\label{eq:di}\,,
\end{align}
with the average $m$ and  $\gamma$ and the maximum deviation amplitudes $\delta m$ and $\delta \gamma$ 
of inertia and damping ratio. Inhomogeneities are determined by the coefficients 
$-1 \le a_i, r_i \le 1$ with $\sum_i r_i=\sum_i a_i=0$ which are determined following a minimization of the performance
measure $\mathcal{M}_b$ of~\eqref{eq:measure}. In the following two paragraphs we construct a matrix perturbation theory
to linear order in the inhomogeneity parameters $\delta m$, and $\delta\gamma$ to calculate
the performance measure $\mathcal{M}_b= \mathcal{M}_b^{(0)} + \sum_ir_i \rho_i +\sum_ia_i \alpha_i + 
\mathcal{O}(\delta m^2, \delta \gamma^2)$. This requires to 
calculate the susceptibilities $\rho_i \equiv \partial \mathcal{M}_b/\partial r_i$
and $\alpha_i \equiv \partial \mathcal{M}_b/\partial a_i$.

\subsection{Inhomogeneity in inertia}

When inertia is inhomogeneous, but the damping ratios remain homogeneous, the system dynamics and  $\mathcal{M}_b$ are still 
given by \eqref{eq:H} and \eqref{eq:Mb}. However, the eigenvectors of the inertia-weighted Laplacian matrix
$\bm{L_M}$ differ from those of $\bm L$ and consequently $\mathcal{M}_b$ is no longer equal to $\mathcal{M}_b^{(0)}$. In general 
there is no simple way to diagonalize $\bm{L_M}$, but one expects that if the inhomogeneity is weak, then the eigenvalues and 
eigenvectors of $\bm{L_M}$ only slightly differ from those of $m^{-1}\bm L$, which allows to construct a perturbation theory.
 \begin{assumption}[Weak inhomogeneity in inertia] The deviations $\delta m \, r_i$ of the local inertias
$m_i$ are all small compared to their average $m$. We write
$\bm M = m\big[\mathbb{1}+\mu\,{\rm diag}\big(\{r_i\}\big)\big]$, where $\mu\equiv\delta m/m\ll 1$
is a small, dimensionless parameter.\label{ass:inertia}
\end{assumption}

To linear order in $\mu$, the series expansion of $\bm{L_M}$ reads
\begin{align}
\bm{L_M}&=\bm M^{-1/2}\bm L\bm M^{-1/2}=m^{-1}\big[\bm L+\mu\bm V_1+\mathcal{O}(\mu^2)\big]\, ,
\end{align}
 with $\bm V_1=-\big(\bm R\bm L+\bm L\bm R\big)/2$ and $\bm R={\rm diag}(\{r_i\})$.
 In this form, the inertia-weighted Laplacian matrix $\bm{L_{M}}$ is given by the sum of 
 an easily  diagonalizable matrix, $m^{-1}\bm L$, and a small 
 perturbation matrix, $(\mu/m) \, \bm V_1$. Truncating the expansion of $\bm{L_M}$ at this linear order 
 gives an error of order $\sim \mu^2$, which is small under Assumption~\ref{ass:inertia}.
 
Matrix perturbation theory gives approximate expressions for the eigenvectors $\bm u_{\alpha}$ 
and eigenvalues $\lambda_\alpha$ of $\bm{L_M}$ in terms of 
those ($\bm u_{\alpha}^{(0)}$ and $\lambda_\alpha ^{(0)}$)   
of $\bm{L}$~\cite{stewart1990matrix}. To leading order in $\mu$ one has
\begin{align}
\lambda_\alpha &= m^{-1}\big[\lambda_\alpha ^{(0)}+\mu\lambda_\alpha ^{(1)}+\mathcal{O}(\mu^2)\big]\,,\label{eq:lambda}\\
\bm u_{\alpha}&= \bm u_{\alpha}^{(0)}+\mu \bm u_{\alpha}^{(1)}+\mathcal{O}(\mu^2)\label{eq:u}\,,
\end{align}
with
\begin{align}
\lambda_\alpha ^{(1)}&=\bm u_{\alpha}^{(0)\top}\bm V_1\bm u_{\alpha}^{(0)}\,,\label{eq:lambda1}\\
\bm u_{\alpha}^{(1)}&=\sum_{\beta\neq \alpha}\frac{\bm u_{\beta}^{(0)\top}\bm V_1 \bm u_{\alpha}^{(0)}}{\lambda_\alpha ^{(0)}-\lambda_\beta^{(0)}}\bm u_{\beta}^{(0)}\,.\label{eq:u1}
\end{align}
From \eqref{eq:Mb}, \eqref{eq:lambda} and \eqref{eq:u}, the first-order approximation of $\mathcal{M}_b$ in $\mu$  reads
\begin{eqnarray}
\mathcal{M}_{b}&=&\mathcal{M}_{b}^{(0)}+\frac{\mu\delta P^2}{2\gamma}\sum_{\alpha>1}\lambda_\alpha ^{(0)-1}\Big(2u_{\alpha b}^{(0)}u_{\alpha b}^{(1)}-r_bu_{\alpha b}^{(0)2}\nonumber\\
&&-u_{\alpha b}^{(0)2}\lambda_\alpha ^{(0)-1}\lambda_\alpha ^{(1)}\Big)+\mathcal{O}(\mu^2)\,.\label{eq:M1}
\end{eqnarray}

\begin{proposition}\label{prop:rho}
For an abrupt power loss,
$\bm{\delta P}(t) = \bm{\delta P} \, \Theta(t)$ on a single bus labeled $b$, 
$\delta P_i=\delta_{ib} \, \delta P$, and
under Assumption~\ref{ass:inertia}, the susceptibilites $\rho_i \equiv \partial \mathcal{M}_b/\partial r_i$  
are given by
\begin{equation}
\rho_i=-\frac{\mu\delta P^2}{\gamma N}\sum_{\alpha>1}\frac{u_{\alpha b}^{(0)}u_{\alpha i}^{(0)}}{\lambda_\alpha ^{(0)}}\, .\label{eq:rho}
\end{equation}
\end{proposition}

\begin{IEEEproof}
Taking the derivative of \eqref{eq:M1} with respect to $r_i$, with $\lambda_\alpha ^{(1)}$ and $u_{\alpha b}^{(1)}$ given in 
\eqref{eq:lambda1} and \eqref{eq:u1}, one gets
\begin{align}
\frac{\partial\mathcal{M}_b}{\partial r_i}&=\frac{\mu\delta P^2}{2\gamma}\bigg[\sum_{\substack{\alpha>1,\\\beta\neq \alpha}}\!u_{\alpha b}^{(0)}u_{\beta b}^{(0)}u_{\alpha i}^{(0)}u_{\beta i}^{(0)}\bigg(\frac{1}{\lambda_\alpha ^{(0)}}-\frac{2}{\lambda_\alpha ^{(0)}-\lambda_\beta^{(0)}}\bigg)\nonumber\\
&-\delta_{ib}\sum_{\alpha>1}\frac{u_{\alpha b}^{(0)2}}{\lambda_\alpha ^{(0)}}+\sum_{\alpha>1}\frac{u_{\alpha b}^{(0)2}u_{\alpha i}^{(0)2}}{\lambda_\alpha ^{(0)}}\bigg]+\mathcal{O}(\mu^2)\,,\label{eq:dmdr}
\end{align}
The first term in the square bracket in \eqref{eq:dmdr} gives
\begin{align}
\!\!\!\sum_{\substack{\alpha>1\\\beta\neq \alpha}}\frac{u_{\alpha b}^{(0)}u_{\beta b}^{(0)}u_{\alpha i}^{(0)}u_{\beta i}^{(0)}}{\lambda_\alpha ^{(0)}} &= \sum_{\substack{\alpha>1,\\\beta}}\frac{u_{\alpha b}^{(0)}u_{\beta b}^{(0)}u_{\alpha i}^{(0)}u_{\beta i}^{(0)}}{\lambda_\alpha ^{(0)}} \nonumber \\
 - \sum_{\alpha>1}\frac{u_{\alpha b}^{(0)2} u_{\alpha i}^{(0)2}}{\lambda_\alpha ^{(0)}}&= \delta_{ib} \sum_{\alpha>1}\frac{u_{\alpha b}^{(0)2}}{\lambda_\alpha ^{(0)}} - \sum_{\alpha>1}\frac{u_{\alpha b}^{(0)2} u_{\alpha i}^{(0)2}}{\lambda_\alpha ^{(0)}} \, , 
\end{align}
where we used $\sum_\beta u_{\beta i}^{(0)}u_{\beta b}^{(0)}=\delta_{ib}$. This terms therefore 
exactly cancels out with the last two terms in the square bracket in \eqref{eq:dmdr} and one obtains 
\begin{equation}
\rho_i(b) = \frac{\partial \mathcal{M}_b}{\partial r_i}=-\frac{\mu\delta P^2}{\gamma}\sum_{\substack{\alpha>1,\\\beta\neq \alpha}}\frac{u_{\alpha b}^{(0)}u_{\beta b}^{(0)}u_{\alpha i}^{(0)}u_{\beta i}^{(0)}}{\lambda_\alpha ^{(0)}-\lambda_\beta^{(0)}}+\mathcal{O}(\mu^2)\,.\label{eq:wk2}
\end{equation}
The argument of the double sum in \eqref{eq:wk2} is odd under permutation of $\alpha$ and $\beta$, therefore only terms with 
$\beta=1$ survive. With $u_{1i}^{(0)}=1/\sqrt{N}$, one finally obtains \eqref{eq:rho}.
\end{IEEEproof}
\begin{remark}\label{rmk:rho}
By summing over every fault locations $b$, one gets $\sum_{b}\rho_i(b)=0$. This follows from the properties of the eigenvector $\bm u_{\alpha}^{(0)}$.
\end{remark}

\subsection{Inhomogeneity in damping ratios}
Equation \eqref{eq:dchii} gives exact solutions to the linearized dynamical problem defined in \eqref{eq:xi_ode}, under the assumption of homogeneous damping ratio, $m_i/d_i \equiv \gamma$. In this section we lift that constraint and write
	$\gamma_i = \gamma + \delta \gamma \, a_i$. With inhomogeneous damping ratios, \eqref{eq:H} becomes
\begin{equation}
\frac{{\rm d}}{{\rm d}t}\left[\!\begin{array}{c}
\bm\xi\\
\bm{\dot\xi}
\end{array}\!\right]=\underbrace{
\left[\!\begin{array}{cc}
\mathbb{0}_{N\times N}&\!\!\!\!\mathbb{1}\\
\!\!\!\!-\bm\Lambda&\!\!\!\!\!\!-\gamma \, \mathbb{1} - \delta \gamma \, \bm V_2
\end{array}\!\right]}_{\bm H}\!\left[\!\begin{array}{c}
\bm\xi\\
\bm{\dot\xi}
\end{array}\!\right]+
\bigg[\!\begin{array}{c}
\mathbb{0}_{N\times 1} \\
\bm{\mathcal{P}}
\end{array}\!\bigg]\,,\label{eq:H1}
\end{equation}
which differs from \eqref{eq:H} only through the additional term $-\delta \gamma \bm V_2$ with
$\bm V_2 = \bm U\bm A\,\bm U^\top$, $\bm A={\rm diag}(\{a_i\})$. Under the assumption that 
that the dimensionless parameter  $g\equiv \delta \gamma/\gamma \ll 1$, 
this additional term gives only small corrections to the unperturbed problem of \eqref{eq:H}, and we use matrix perturbation 
theory to calculate these corrections in a polynomial expansion in $g$.

\begin{assumption}[Weak inhomogeneity in damping ratios] The deviations $\delta \gamma \, a_i$ of the damping ratio
$\gamma_i$ from their average $\gamma$ are all small compared to their average. We write
$\bm \Gamma = \gamma\big[\mathbb{1}+g\,{\rm diag}\big(\{a_i\}\big)\big]$, where $g\equiv\delta\gamma/\gamma \ll 1$
is a small, dimensionless parameter.\label{ass:gamma}
\end{assumption}
We want to integrate \eqref{eq:H1} using the spectral approach that provided the solutions 
\eqref{eq:chiipm}. In principle this requires to know the eigenvalues and eigenvectors of $\bm H$ in \eqref{eq:H1}, 
which is not possible in general, because $\bm V_2$ does not commute with $\bm \Lambda$. 
When $g$ is small enough, the eigenvalues and eigenvectors are only slightly 
altered~\cite{stewart1990matrix} and can be systematically calculated order
by order in a polynomial expansion in $g$. We therefore follow a perturbative approach which expresses solutions
to \eqref{eq:H1} in such a polynomial expansion in $g$. Formally, one has, for the eigenvalues $\mu_{\alpha s}$
and for the left and right eigenvectors $\bm t_{\alpha s}^{L,R}$ of $\bm H$
\begin{align}
\mu_{\alpha s}&=\sum_{m=1}^\infty g^m \, \mu_{\alpha s}^{(m)}  \label{eq:mu_exp} \, , \\
\bm{t}_{\alpha s }^{L,R}&= \sum_{m=1}^\infty g^m 
\, \bm{t}_{\alpha s}^{(m)L,R} \label{eq:TL_exp}\, ,
\end{align}
where the $m=0$ terms are given by the eigenvalues, $\mu_{\alpha s}^{(0)}$, and the left and right eigenvectors, 
$\bm{t}_{\alpha s}^{(0)L,R}$, of the matrix $\bm H_0$  in \eqref{eq:H}, corresponding to homogeneous inertia.
In order for the sums in \eqref{eq:mu_exp} and \eqref{eq:TL_exp} to converge, a necessary condition is that $g <1$. 
The task is to calculate the terms $\mu_{\alpha s}^{(m)}$ and $\bm{t}_{\alpha s}^{(m)L,R}$ with $m=1,2,...$.
When $g \ll 1$, one expects that only few, low order terms already give a good estimate 
of the eigenvalues and eigenvectors of $\bm H$. In this manuscript, we calculate the first-order corrections, $m=1$. 
They are given by formulas similar to \eqref{eq:lambda1} and \eqref{eq:u1},
\begin{align}
& g \, \mu_{\alpha s}^{(1)}= \bm{t}_{\alpha s}^{(0)L}  \left[\!\begin{array}{cc}
\mathbb{0}_{N\times N}&\mathbb{0}_{N\times N}\\
\mathbb{0}_{N\times N}&- \delta \gamma \, \bm V_2
\end{array}\!\right]  \bm{t}_{\alpha s}^{(0)R} \, , \\
&g \, \bm{t}_{\alpha s}^{(1)R}=\overline{\sum_{\beta,s'}} \; \frac{\bm{t}_{\beta s'}^{(0)L}  \left[\!\begin{array}{cc}
\mathbb{0}_{N\times N}&\mathbb{0}_{N\times N}\\
\mathbb{0}_{N\times N}&- \delta \gamma \, \bm V_2
\end{array}\!\right]  \bm{t}_{\alpha s}^{(0) R}}{\mu_{\alpha s}^{(0)}-\mu_{\beta s'}^{(0)}} \bm{t}_{\beta s'}^{(0) R} \, , \\
&g \, \bm{t}_{\alpha s}^{(1)L}=\overline{\sum_{\beta,s'}} \; \frac{\bm{t}_{\alpha s}^{(0)L}  \left[\!\begin{array}{cc}
\mathbb{0}_{N\times N}&\mathbb{0}_{N\times N}\\
\mathbb{0}_{N\times N}&- \delta \gamma \, \bm V_2
\end{array}\!\right]  \bm{t}_{\beta s'}^{(0) R}}{\mu_{\alpha s}^{(0)}-\mu_{\beta s'}^{(0)}} \bm{t}_{\beta s'}^{(0) L} \, ,
\end{align}
where 
${\overline\sum}$ indicates that the sum runs over $(\beta,s') \ne (\alpha, s)$. 
One obtains
\begin{align}
&g \, \mu_{\alpha s}^{(1)} = -\delta \gamma\Big(\frac{1}{2}+ i s \frac{\gamma}{2 f_\alpha }\Big)\bm{V}_{2;\alpha\alpha}\,, \!\!\label{eq:pertval} \\
&g \, \bm{t}_{\alpha s}^{(1)R}=  2 \, \delta \gamma \, \overline{\sum_{\beta,s'}} \frac{\bm{V}_{2;\alpha\beta} \, \mu_{\alpha s}^{(0)}}{f_\beta (s s' \, f_\alpha-f_\beta)} \, \bm{t}_{\beta s'}^{(0) R} \, , 
\label{eq:pertvec1}\\
&g \, \bm{t}_{\alpha s}^{(1)L}= 2 \, \delta \gamma \, \overline{\sum_{\beta,s'}} \frac{\bm{V}_{2;\alpha\beta} \, \mu_{\beta s'}^{(0)}}{f_\alpha (f_\alpha-s s' \, f_\beta)} \, \bm{t}_{\beta s'}^{(0) L}  \, ,\label{eq:pertvec}
\end{align}
with $\bm{V}_{2;\alpha\beta}=\sum_{i}a_i\,u_{\alpha i} \,u_{\beta i}$.

\begin{remark}
By definition, $-1 \le \bm{V}_{2;\alpha\alpha} \le 1$. Therefore,
\eqref{eq:pertval} indicates, among others,
that when the parameters $\{a_i\}$ are correlated (anticorrelated) with the square components 
$\{u_{\alpha i}^2\}$ for some $\alpha$ then that mode is more strongly (more weakly) damped.
Accordingly, Theorem~\ref{thm:lin_opt_control} will distribute the set $\{a_i\}$ to increase the damping of the 
slow modes of $\bm H$. 
\end{remark}

\begin{proposition}
For an abrupt power loss,
$\bm{\delta P}(t) = \bm{\delta P} \, \Theta(t)$ on a single bus labeled $b$, 
$\delta P_i=\delta_{ib} \, \delta P$,
and under Assumption \ref{ass:gamma}, $\dot\xi_\alpha(t)$ reads, to leading order
in $g$,
\begin{align}
\dot\xi_\alpha(t)=&\frac{\mathcal{P}_\alpha }{f_\alpha }e^{-\gamma t/2}\bigg[2s_\alpha \Big( 1+g \frac{\gamma^2}{f_\alpha^2} \bm{V}_{2; \alpha \alpha} \Big)\nonumber\\
& \;\;\;\;\;\;\;\;\;\;\;\;\;\;\;\;\;\;-g\gamma t\bm{V}_{2;\alpha\alpha}\Big(s_\alpha+\frac{\gamma}{f_\alpha }c_\alpha\Big)\bigg]\nonumber\\
&+g\gamma\sum_{\beta\neq \alpha}\frac{\bm{V}_{2;\alpha\beta}\mathcal{P}_\beta}{\lambda_\alpha -\lambda_\beta}e^{-\gamma t/2}\bigg[\frac{\gamma}{f_\beta}s_\beta-\frac{\gamma}{f_\alpha }s_\alpha+c_\alpha-c_\beta\bigg]\nonumber\\
&+\mathcal{O}(g^2)\,,\label{eq:perturbed_xii}
\end{align}\label{prop:perturbed_xii}
where $s_\alpha=\sin(f_\alpha t/2)$ and $c_\alpha=\cos(f_\alpha t/2)$, and $\mathcal{P}_\alpha$ and $f_\alpha$ are defined below
\eqref{eq:dchii}.
\end{proposition}
The proof is based on \eqref{eq:pertval} to \eqref{eq:pertvec} and is given in Appendix~\ref{sect:continuation}.

\begin{proposition}\label{prop:alpha}
For an abrupt power loss,
$\bm{\delta P}(t) = \bm{\delta P} \, \Theta(t)$ on a single bus labeled $b$, 
$\delta P_i=\delta_{ib} \, \delta P$, and
under Assumption~\ref{ass:gamma}, the susceptibilities $\alpha_i \equiv \partial \mathcal{M}_b/\partial a_i$ are given by
\begin{eqnarray}\label{eq:alphai}
\alpha_i &=& - \frac{g\delta P^2}{2\gamma m_b} \Bigg[ \sum_{\alpha>1} \frac{u_{\alpha i}^{2}u_{\alpha b}^{2}}{\lambda_\alpha } 
 \nonumber \\
&& +  \sum_{\substack{\alpha>1,\\\beta \ne \alpha}} \frac{
\, u_{\alpha i}\, u_{\alpha b} u_{\beta i}\, u_{\beta b} }
{(\lambda_\alpha-\lambda_\beta)^2 + 2 \gamma^2 (\lambda_\alpha+\lambda_\beta) }\Bigg] \, \qquad
\end{eqnarray}
\end{proposition}
\begin{IEEEproof}
From \eqref{eq:perturbed_xii}, to first order in $g$, one has
\begin{align}
&\intop_0^{\infty}\dot\xi_\alpha^2(t){\rm d}t=\frac{\mathcal{P}_\alpha ^2}{2\gamma\lambda_\alpha }\left(1-g \bm{V}_{2;\alpha\alpha}  \right)\nonumber\\
&-g\gamma \sum_{\beta\neq \alpha}\frac{\bm{V}_{2;\alpha\beta} \, \mathcal{P}_\alpha\mathcal{P}_\beta}{(\lambda_\alpha -\lambda_\beta)^2+2\gamma^2(\lambda_\alpha +\lambda_\beta)} 
+\mathcal{O}(g^2)\,.\label{eq:intop_xi}
\end{align}
Taking the derivative of \eqref{eq:intop_xi} with respect to $a_i$ with the definition of $\bm{V}_{2;\alpha\beta}$ given below 
\eqref{eq:pertvec}, and summing over $\alpha >1$ one obtains \eqref{eq:alphai}.
\end{IEEEproof}
\begin{remark}
We have found numerically that the second term is generally much smaller than the first one and gives only marginal corrections
to our optimized solution. 
\end{remark}
\begin{remark}\label{rmk:alpha}
Close to the homogeneous case $\bm M=m\mathbb{1}$ and $\bm \Gamma=\gamma\mathbb{1}$, summing over every fault locations $b$ makes the second term in \eqref{eq:alphai} vanish. One gets $\sum_b\alpha_i = -g\delta P^2\sum_{\alpha>1}
u_{\alpha i}^{(0)2}/(2\gamma\lambda_\alpha ^{(0)})$. This follows from the properties of the eigenvector of $\bm u_{\alpha}^{(0)}$.
\end{remark}

\section{Optimal placement of inertia and primary control}\label{section:optimal}

In general it is not possible to obtain closed-form 
analytical expressions for the parameters $a_i$ and $r_i$ determining the 
optimal placement of inertia and primary control. Simple optimization algorithms can however be constructed that 
determine how to distribute these parameters to minimize $\mathcal{M}_b$. Theorems~\ref{thm:lin_opt_inertia} and
\ref{thm:lin_opt_control} give two such algorithms for optimization under Assumption~\ref{ass:inertia}  and \ref{ass:gamma} 
respectively. Additionally, Conjecture~\ref{conj:opt} proposes an algorithm for optimization under both 
Assumption~\ref{ass:inertia}  and \ref{ass:gamma}.

\begin{theorem}
\label{thm:lin_opt_inertia}
For an abrupt power loss, under Assumption \ref{ass:inertia} and with $\bm \Gamma=\gamma\mathbb{1}$, the optimal distribution
of parameters $\{r_i\}$ that minimizes $\mathcal{M}_b$ is obtained as follows.
\begin{enumerate}
\item{Compute the sensitivities}
$\rho_i= \partial \mathcal{M}_b/\partial r_i$ from \eqref{eq:rho}
\item{Sort the set $\{\rho_i\}_{i=1, ... N}$ in ascending order}
\item{Set $r_i=1$ for $i=1,... {\rm Int}[N/2]$ and $r_i=-1$ for  $i=N-{\rm Int}[N/2]+1,...N$}
\end{enumerate}
The optimal placement of inertia and primary control is given by
\begin{align}
m_i=m+\delta m\, r_i\,,\;\;\;\;\;\;\;
d_i=\gamma(m+\delta m\, r_i)\,.
\end{align}
\end{theorem}
The proof is in Appendix~\ref{sect:continuation}.

\begin{theorem}
\label{thm:lin_opt_control}
For an abrupt power loss, under Assumption \ref{ass:gamma} and with $\bm M=m\mathbb{1}$, the optimal distribution of parameters 
$\{a_i\}$ that minimizes $\mathcal{M}_b$ is obtained as follows.
\begin{enumerate}
\item{Compute the sensitivities}
$\alpha_i=\partial \mathcal{M}_b/\partial r_i$ from \eqref{eq:alphai}.
\item{Sort the set $\{\alpha_i\}$ in ascending order,}
\item{Set $a_i=1$ for $i=1,... {\rm Int}[N/2]$ and for  $i=N-{\rm Int}[N/2]+1,...N$}
\end{enumerate}
The optimal placement of primary control is given by
\begin{equation}
d_i=m(\gamma+\delta\gamma\, a_i)\,.
\end{equation}
\end{theorem}
\begin{IEEEproof}
With Proposition~\ref{prop:alpha} and $\bm M=m\mathbb{1}$, we get \eqref{eq:alphai}. The proof is the same as the one for Theorem~\ref{thm:lin_opt_inertia} given in Appendix~\ref{sect:continuation}, but with $\{\alpha_i\}$ instead of $\{\rho_i\}$ .
\end{IEEEproof}
We next conjecture an algorithmic combined linear optimization treating simultaneously  Assumptions~\ref{ass:inertia} and 
\ref{ass:gamma}. The difficulty is that for fixed total inertia and damping, one must have $\sum_i m_i=N \, m $, $\sum_i d_i=N \,d$.
From \eqref{eq:di}, the second condition requires $\sum_i a_i r_i = 0$. This is a quadratic, nonconvex constraint, which makes the 
problem nontrivial to solve. The following conjecture presents an algorithm that starts from the distribution 
$\{a_i\}$ and $\{r_i\}$ from Theorems~\ref{thm:lin_opt_inertia} and \ref{thm:lin_opt_control} and orthogonalizes them
while trying to minimize the related increase in $\mathcal{M}_b$.
\begin{conjecture}[Combined linear optimization]
\label{conj:opt} For an abrupt power loss, under Assumptions~\ref{ass:inertia} and \ref{ass:gamma}, the optimal placement of a fixed
 total amount of inertia $\sum_i m_i=mN$ and primary control $\sum_i d_i=dN$ that minimize $\mathcal{M}_b$ is obtained as 
 follows.
\begin{enumerate}
\item{Compute the parameters 
$r_i$ and $a_i$ from Theorems \ref{thm:lin_opt_inertia} and \ref{thm:lin_opt_control}.}
\begin{enumerate}
\item{If $N$ is odd, align the zeros of $\{r_i\}$ and $\{a_i\}$ 
. Let $i_{r0}$ and  $i_{a0}$ be the indexes of these zeros. Their new common index is}
$$
i_{\rm align} = \underset{i}{\rm argmin}(r_i\rho_{i_{r0}}+a_i\alpha_{i_{a0}}-r_i\rho_{i}-a_i\alpha_{i})\,.
$$
Interchange the parameter values $r_{i_{r0}} \leftrightarrow r_{i_{\rm align}}$ and $a_{i_{r0}} \leftrightarrow a_{i_{\rm align}}$.
\item{If $N$ is even, do nothing}
\end{enumerate}
\item{If $n\equiv\sum_ir_ia_i=0$, the optimization is done.}
\item{Find the set $\mathcal{I}=\{i\, | \, {\rm sgn}(r_ia_i)={\rm sgn}(n) \}$. To reach $\sum_ir_ia_i \rightarrow 0$, our strategy is to
set to zero some elements of $\mathcal{I}$. Since however $\sum_i a_i = \sum_i r_i =0$ must be conserved, this must 
be accompanied by a simultaneous change of some other parameter.}
\item{Find the pair $(a_{i1},a_{i2}=-a_{i1})$ or $(r_{i1},r_{i2}=-r_{i1}) \in \mathcal{I} \times \mathcal{I}$ 
which, when sent to $(0,0)$, induce the smallest increase of the objective function $\mathcal{M}_b$. 
Send it to $(0,0)$. Because the pair has opposite sign, this does not affect the condition $\sum_i a_i = \sum_i r_i =0$.}
\item{go to step \# 2.}
\end{enumerate}
\end{conjecture}
It is not at all guaranteed that the algorithm presented in Conjecture~\ref{conj:opt} is optimal, however, numerical results to be
presented below indicate that it works well. 

The optimization considered so far focused on a single fault on bus labeled $b$. We are interested, however, in finding the optimal
distribution of inertia and/or primary control for all possible faults. To that end we introduce the following global 
vulnerability measure 
\begin{equation}
\mathcal{V}=\sum_b\eta_b \, \mathcal{M}_b(\delta P_b)\,,\label{eq:vulnerability}
\end{equation}
where the sum runs over all generator buses. The vulnerability measure $\mathcal{V}$
gives a weighted average over all possible fault positions, with the weight $\eta_b$ 
accounting for the probability that a fault occurs at $b$ and $\delta P_b$ accounting for its potential intensity as given, e.g. by
the rated power of the generator at bus $b$.

\begin{figure}[h!]
\center
\includegraphics[width=0.15\textwidth]{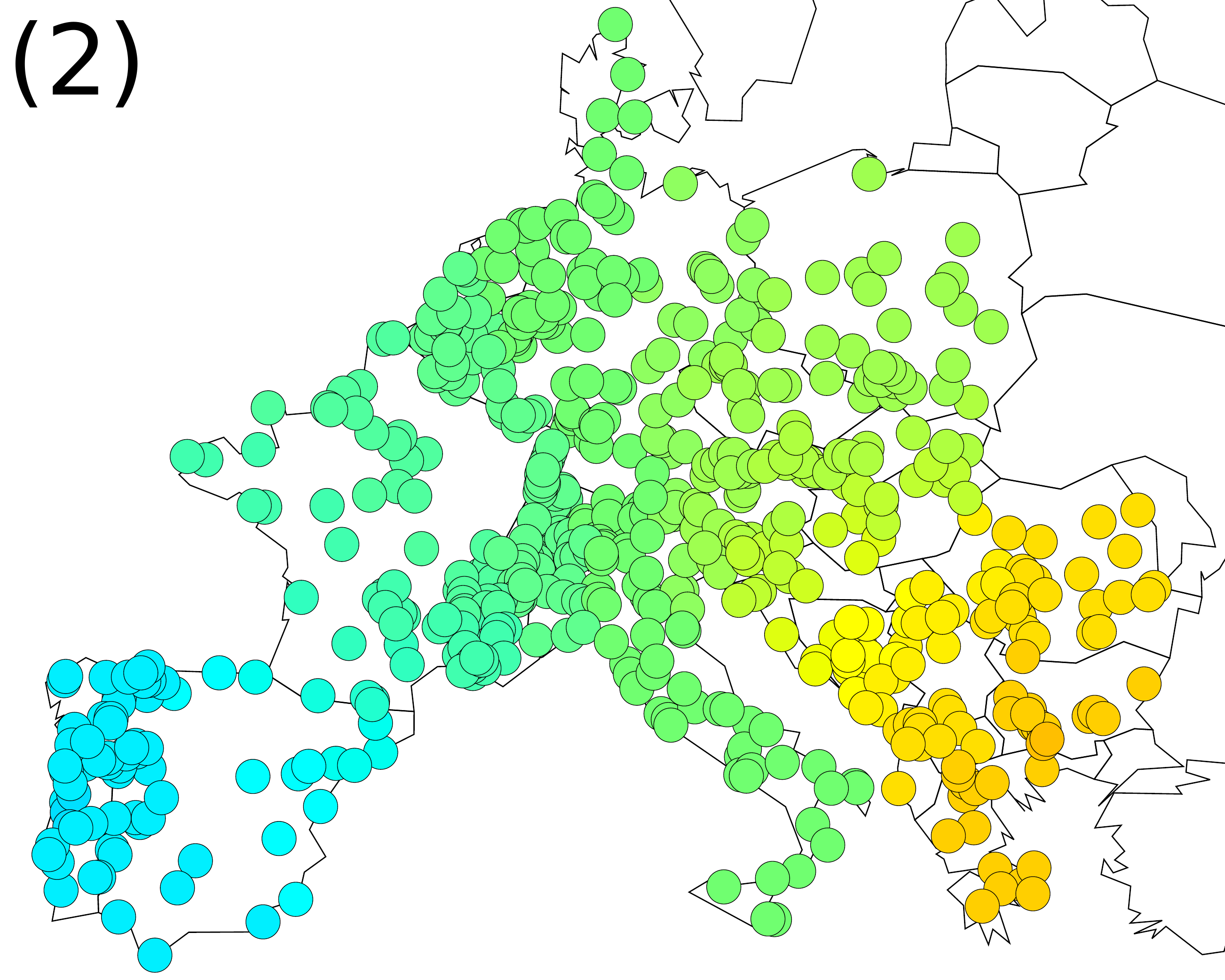}
\includegraphics[width=0.15\textwidth]{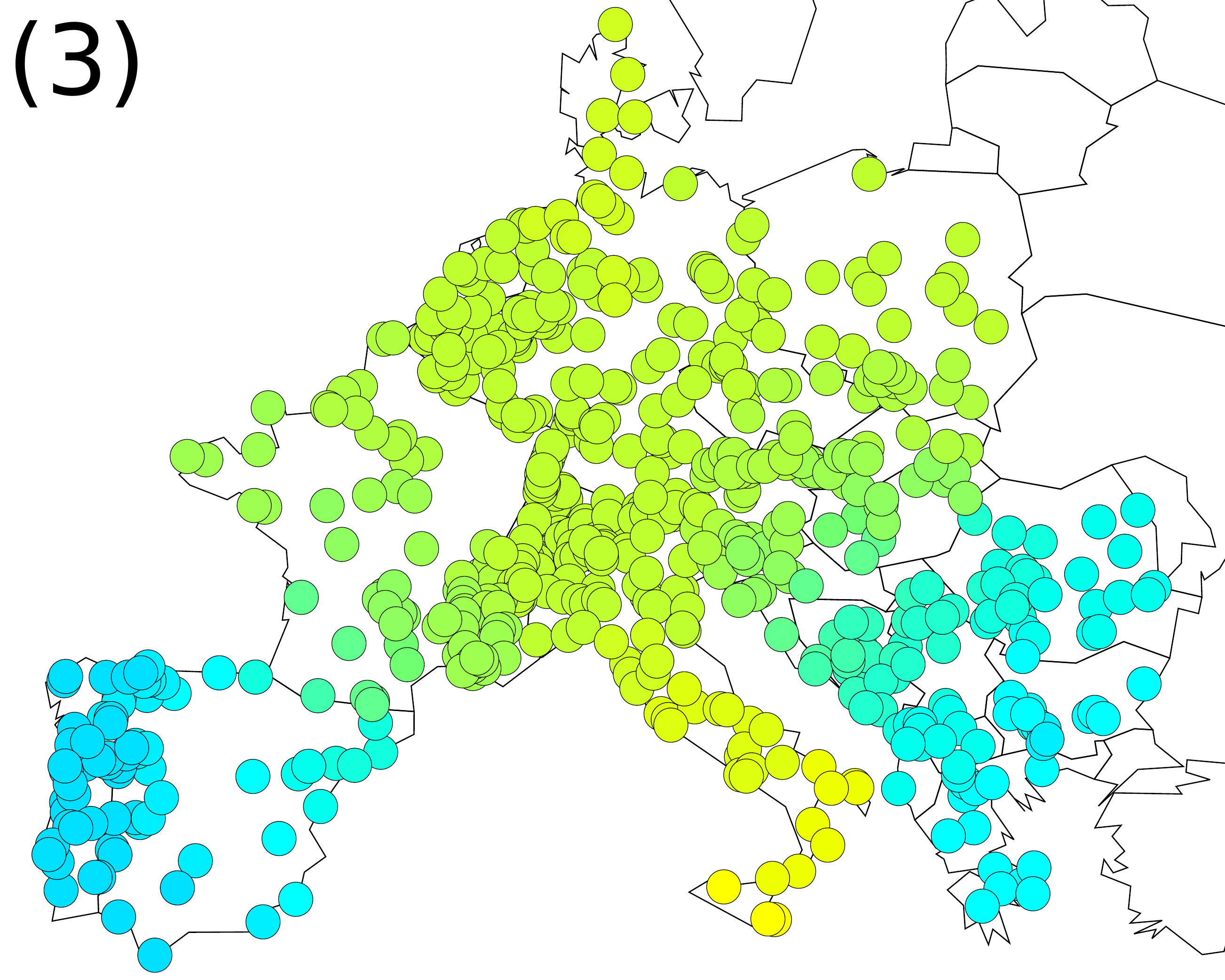}
\includegraphics[width=0.15\textwidth]{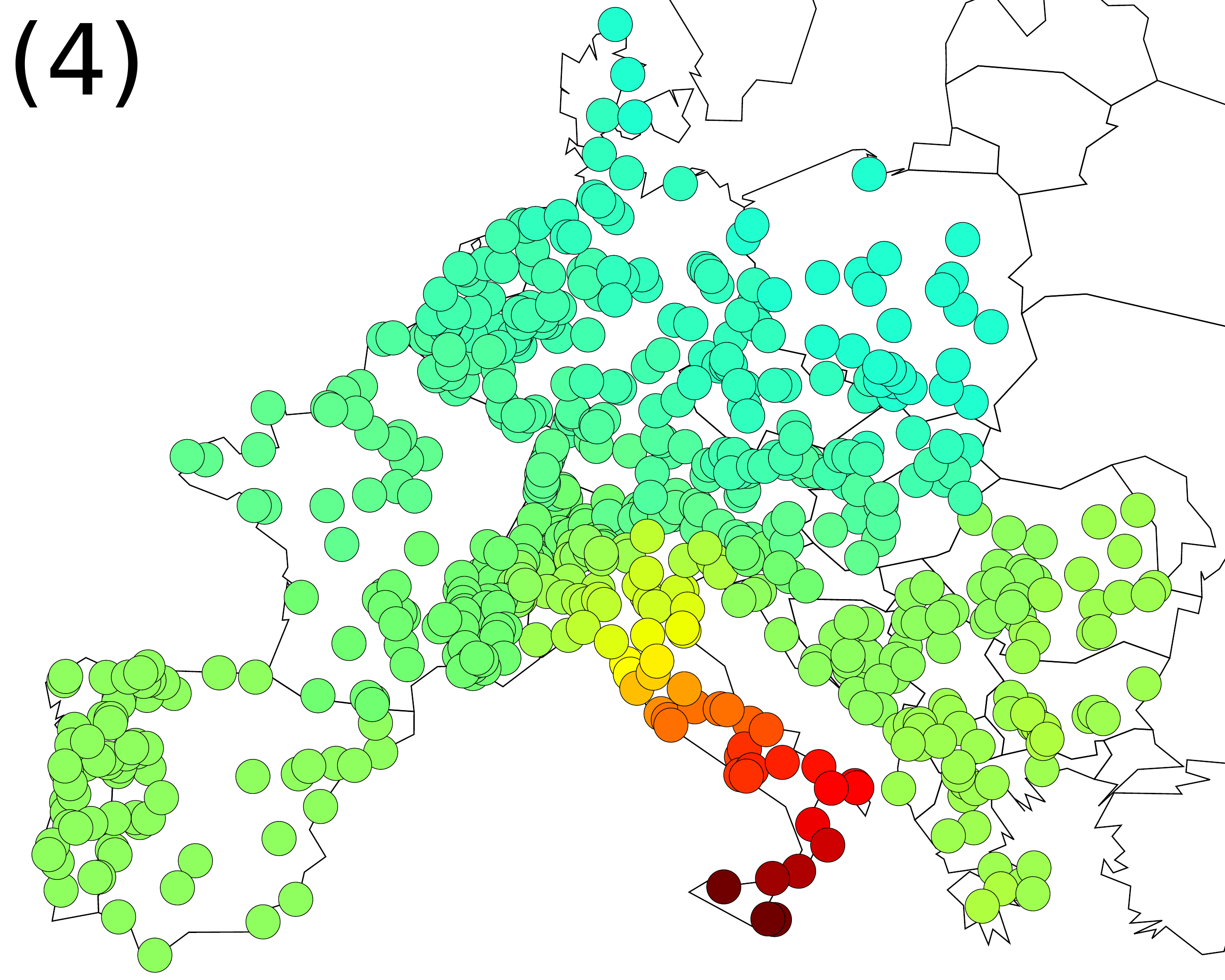}\\
\includegraphics[width=0.15\textwidth]{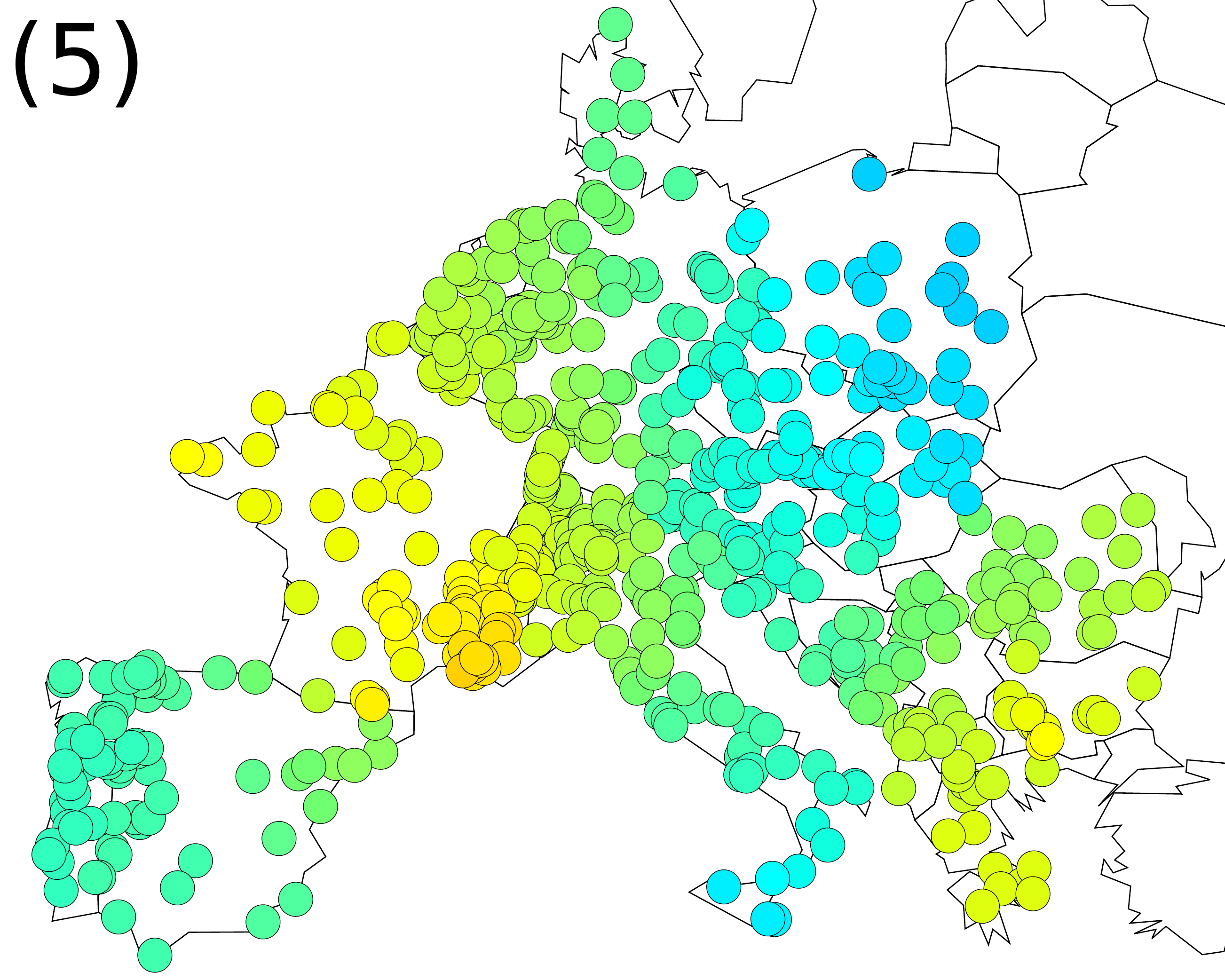}
\includegraphics[width=0.15\textwidth]{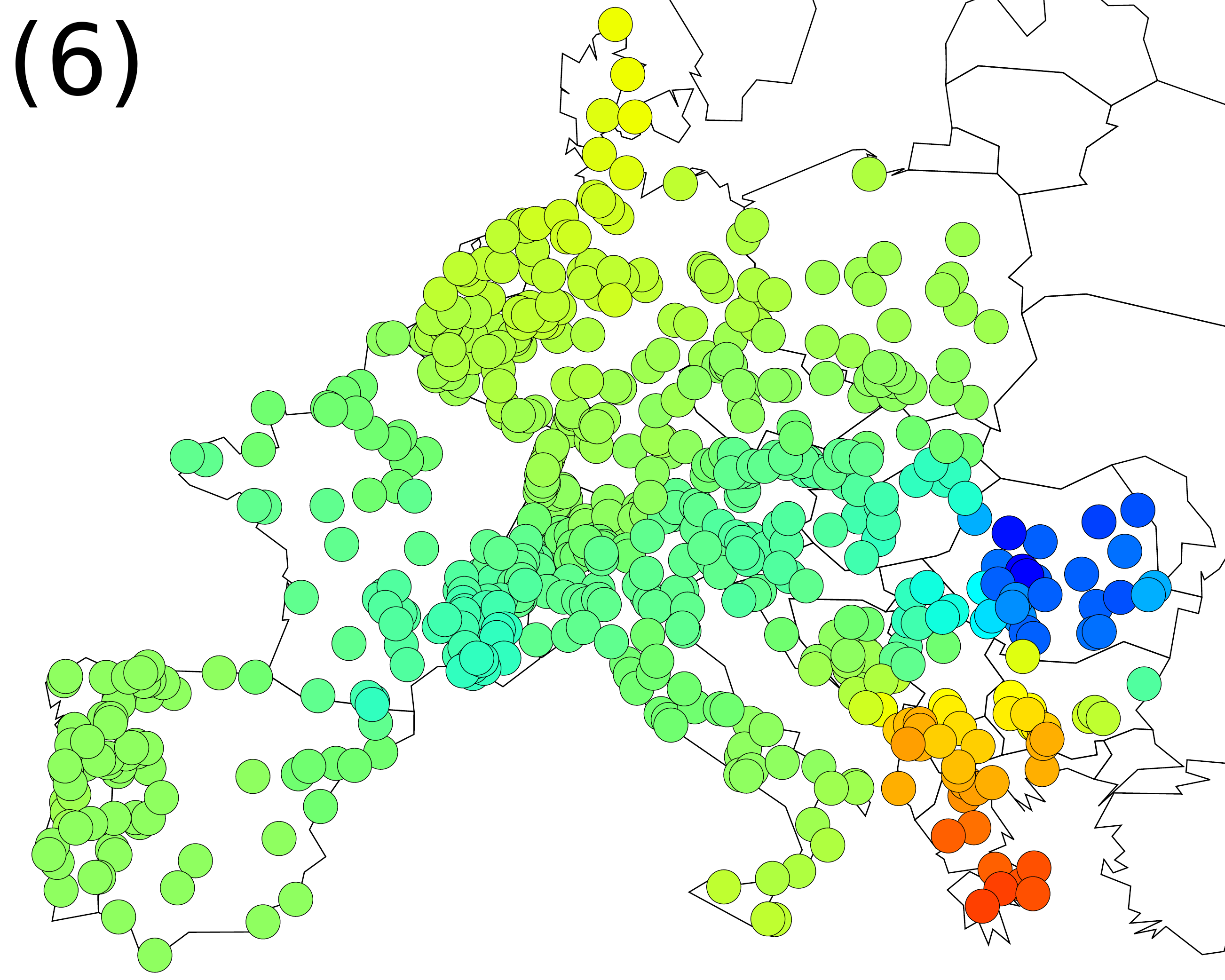}
\includegraphics[width=0.15\textwidth]{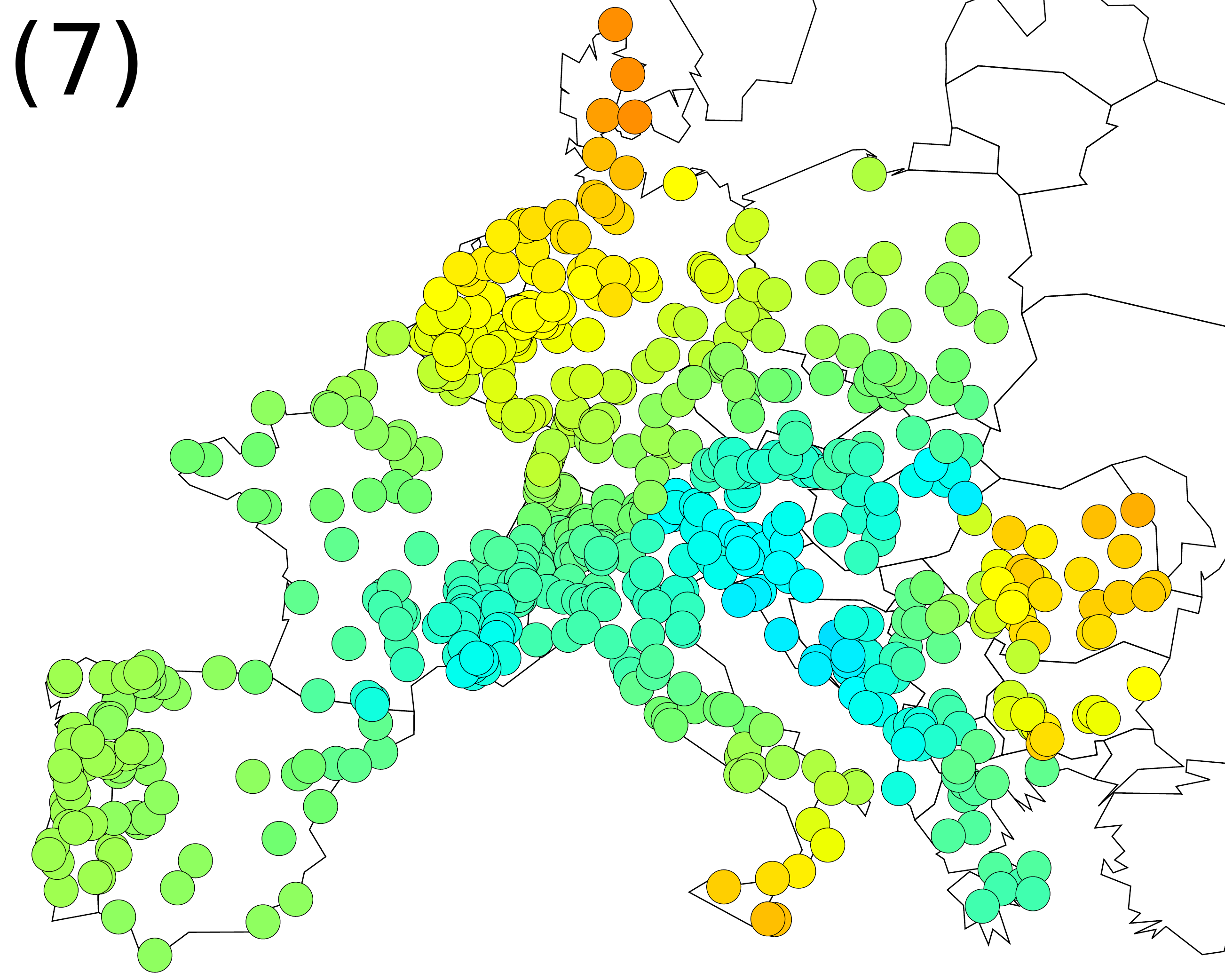}\\
\includegraphics[width=0.27\textwidth]{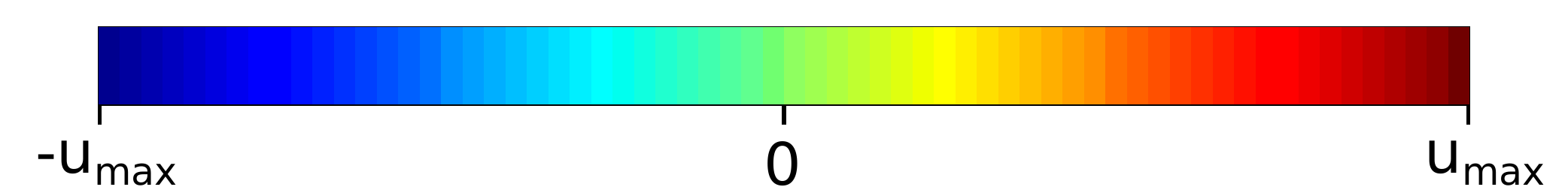}\\
\caption{Color-coded components of the $\alpha=2,3,...7$ eigenvectors of $\bm L$. The colors
span the interval $[-u_{\max},u_{\max}]$ where $u_{\max}=\max_{\alpha\in\{2,\cdots,7\}}\big|u_{\alpha i}^{(0)}\big|$.}\label{fig:slow_modes}
\end{figure}

For equiprobable fault locations and for the same power loss everywhere, $\eta_b\equiv1$, with Remark~\ref{rmk:rho}, it is straightforward to see that 
$\partial \mathcal{V}/\partial r_i = 0+\mathcal{O}(\mu^2)$. 
Therefore, to leading order, there is no benefit in scaling up the inertia anywhere. 
On the other hand, with Remark~\ref{rmk:alpha}, we get $\partial \mathcal{V}/\partial a_i = -g\delta P^2\sum_{\alpha>1}
u_{\alpha i}^{(0)2}/(2\gamma\lambda_\alpha ^{(0)})+\mathcal{O}(g^2)$. The corresponding optimal placement of primary control can 
be obtained with Theorem~\ref{thm:lin_opt_control}, from which we observe that the damping ratios are increased for the buses 
with large squared components $u_{\alpha i}^{(0)2}$ of the slow modes of $\bm L$ -- those with the smallest $\lambda_\alpha^{(0)}$. These modes are displayed in 
Fig~\ref{fig:slow_modes}. One  concludes that, with a non-weighted vulnerability measure, $\eta_b\equiv1$
in \eqref{eq:vulnerability}, an homogeneous inertia location is a local optimum for $\mathcal{V}$, for which damping parameters
need to be increased primarily on peripheral buses. 

\section{Numerical investigations}\label{section:applications}

We illustrate our main results on a model of the synchronous power grid of continental Europe. The network
has 3809 nodes, among them 618 generators, connected through 4944 lines. For details of the model and its
construction we refer the reader to~\cite{pagnier2019inertia,tyloo2018key}. 
To connect to the theory presented above, we remove inertialess buses through 
a Kron reduction \cite{dorfler2013kron} and uniformize the distribution of inertia to $m_i=29.22$MWs$^2$, 
 and primary control $d_i=12.25$MWs. This guarantees that the total amounts of inertia and primary control 
are kept at their initial levels.

\begin{figure}[h!]
\center
\includegraphics[width=0.24\textwidth]{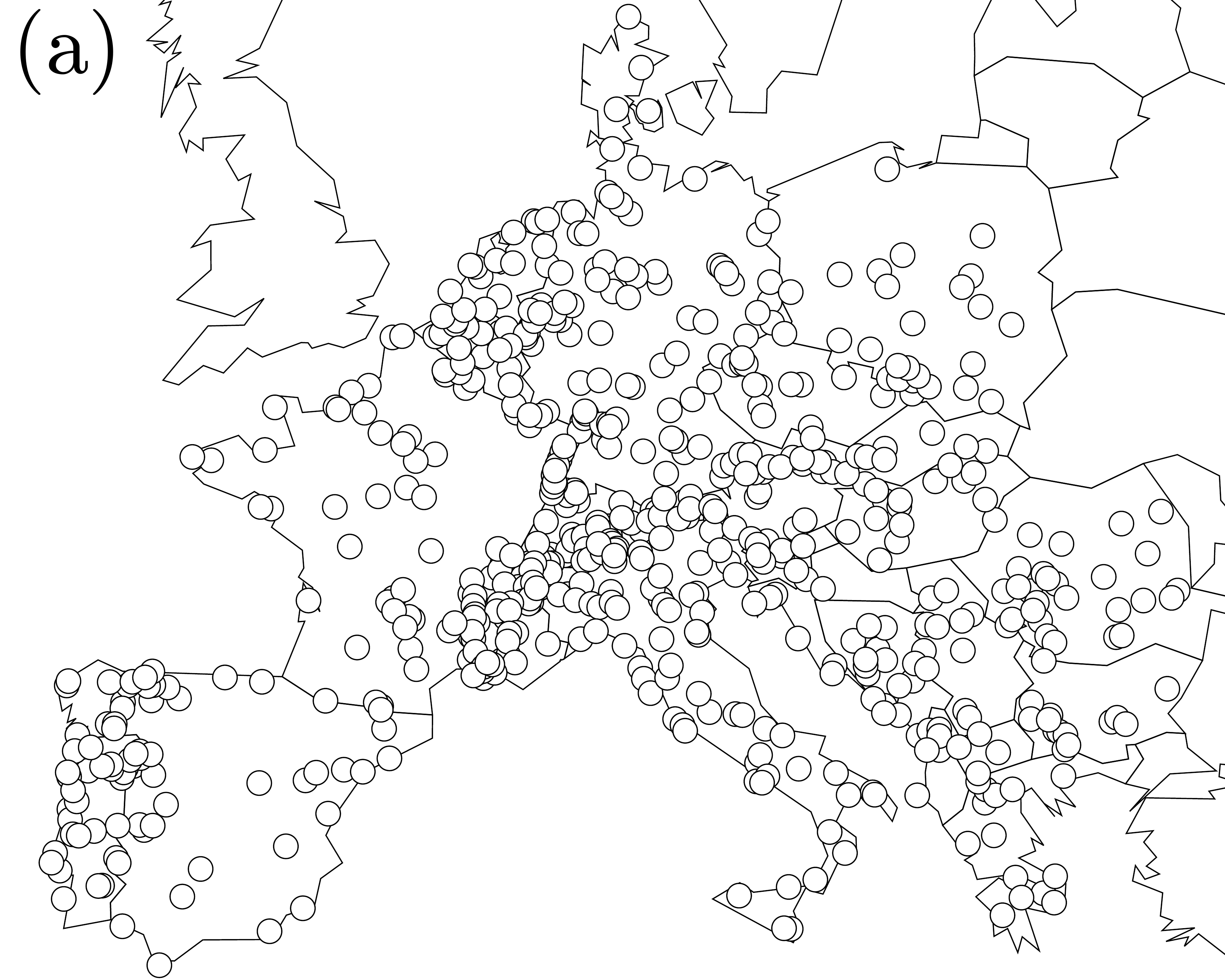}
\includegraphics[width=0.24\textwidth]{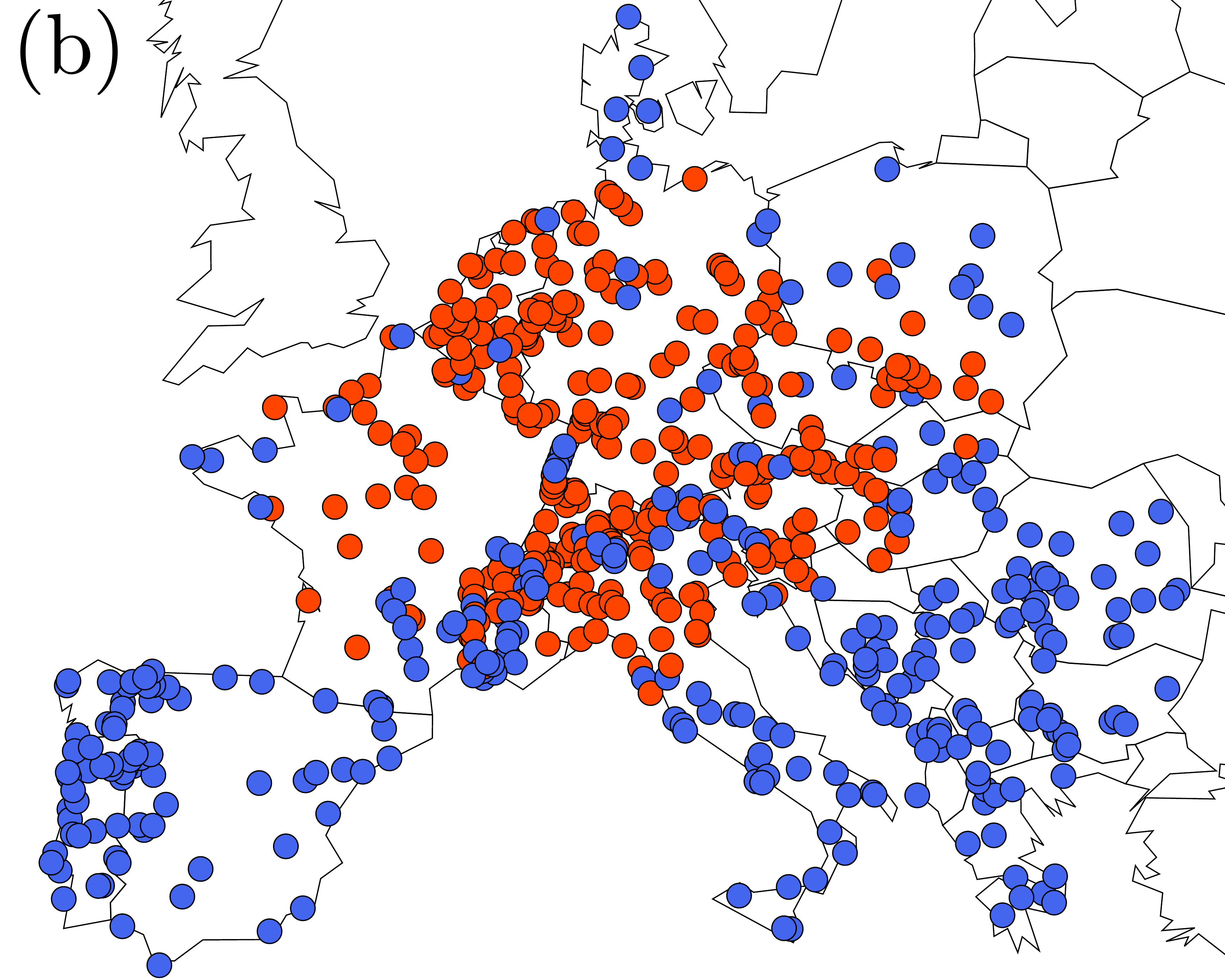}\\
\includegraphics[width=0.24\textwidth]{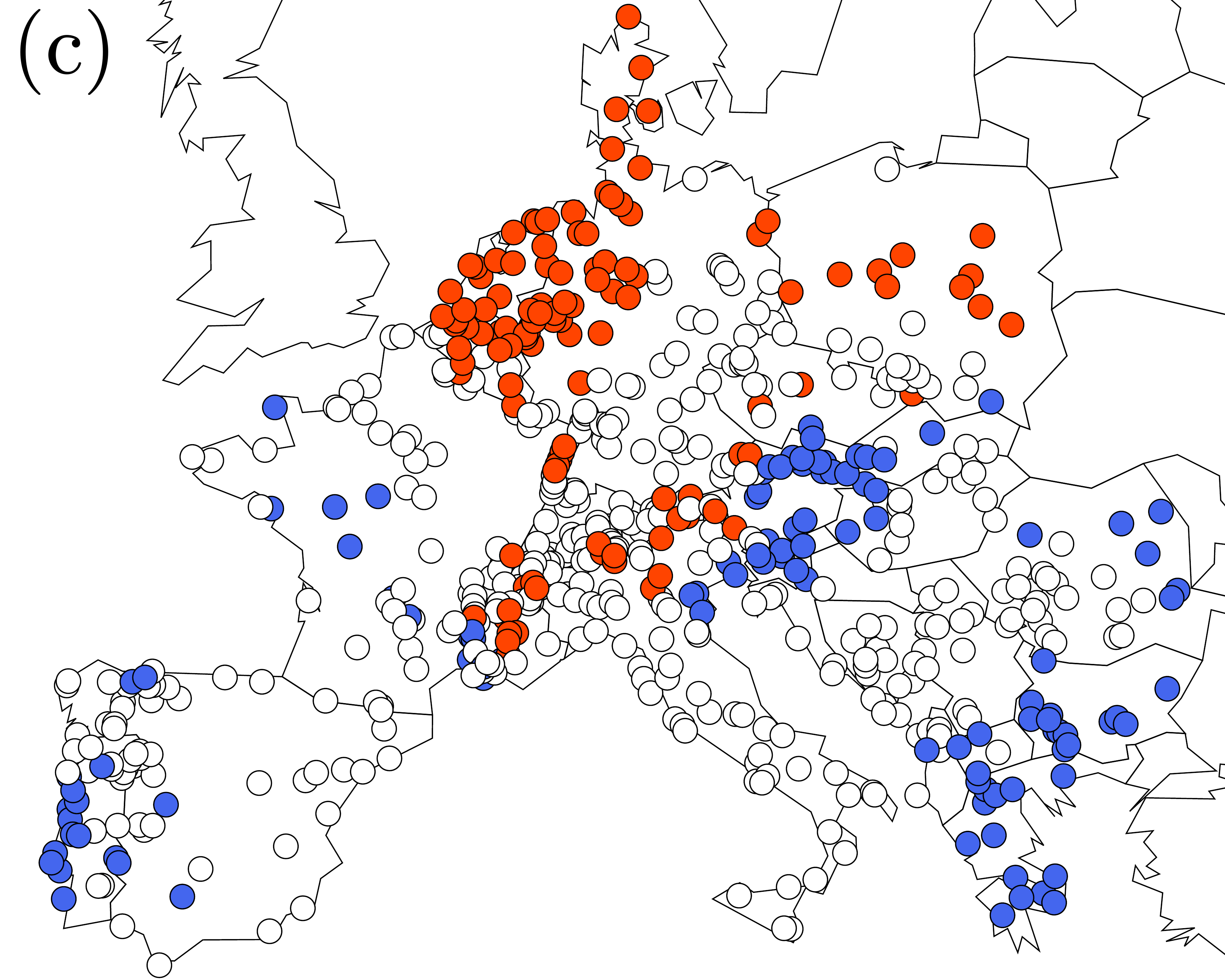}
\includegraphics[width=0.24\textwidth]{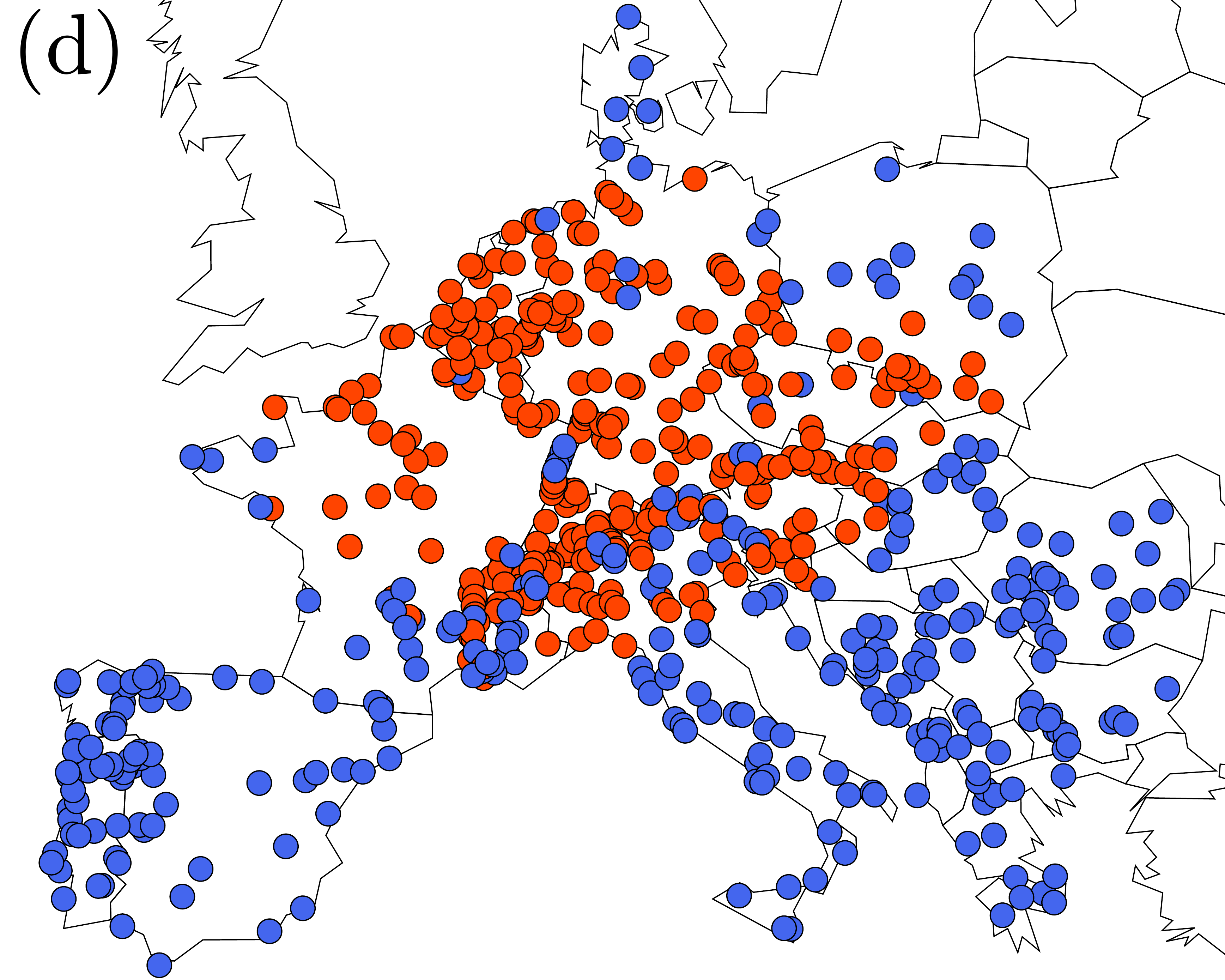}\\
\includegraphics[width=0.24\textwidth]{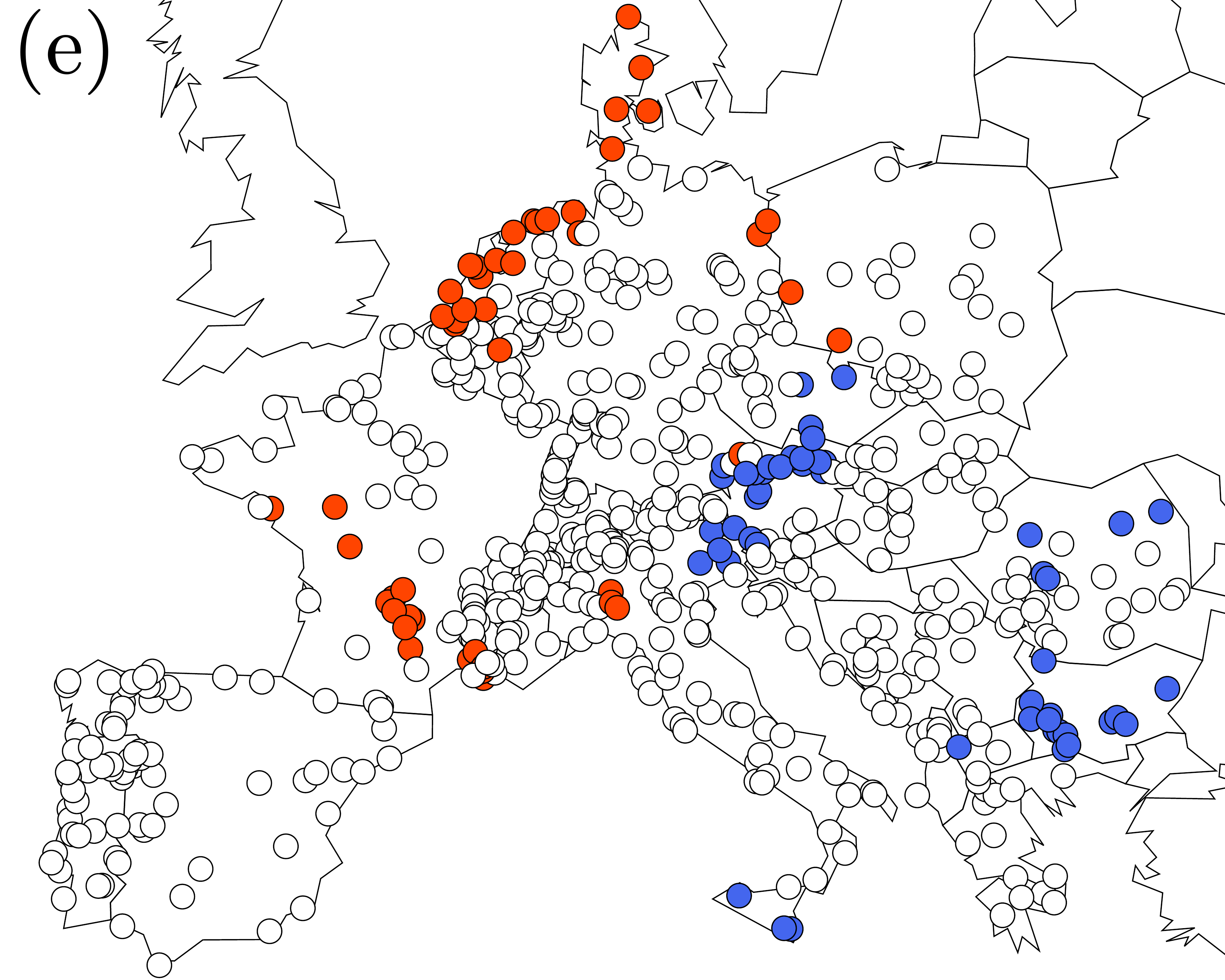}
\includegraphics[width=0.24\textwidth]{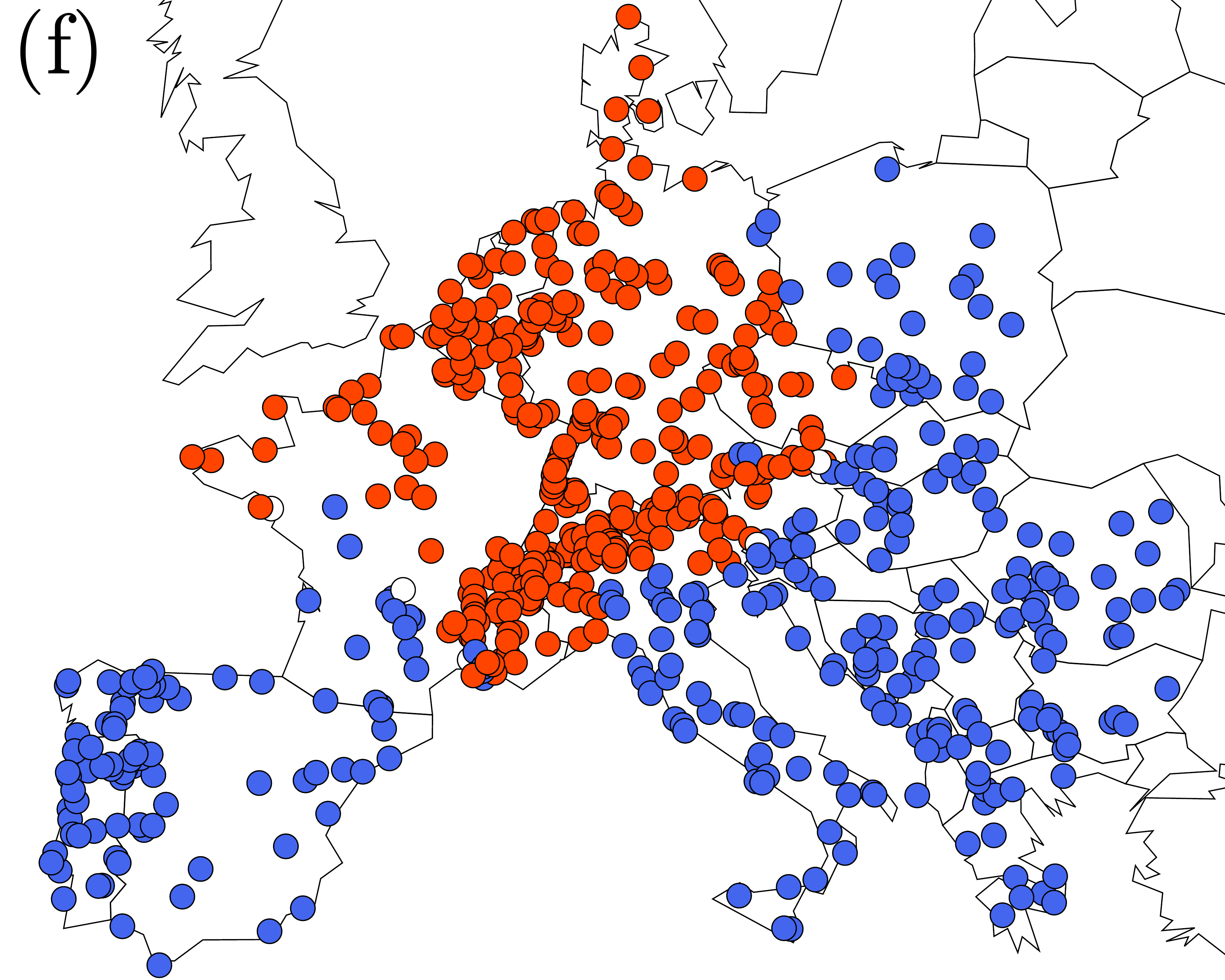}\\
\vspace{5pt}
\includegraphics[width=0.47\textwidth]{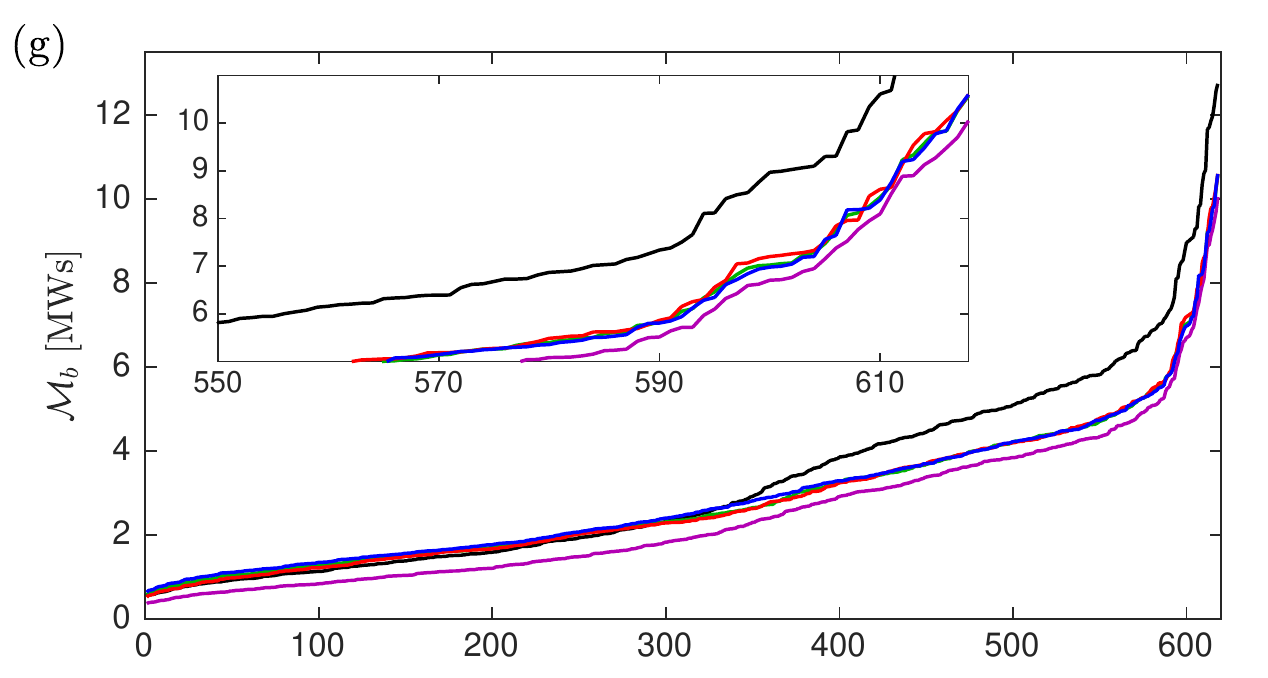}
\caption{Deviation from homogeneous inertia and primary control following the minimization of $\mathcal{V}$ in \eqref{eq:vulnerability} for different choices, (a)-(b) $\eta_b\equiv 1$, (c)-(d) $\eta_b=\mathcal{M}_b^{(0)2}$ and (e)-(f) as defined in \eqref{eq:thres}. $r_i=-1,0,1$ (left) and $a_i=-1,0,1$ (right) are displayed in red, white and blue respectively. (g) Vulnerability $\mathcal{M}_b$ vs. fault location (in increasing order of $\mathcal{M}_b$) for the homogeneous model (black) and the optimized models corresponding to (a)-(b) (green line), (c)-(d) (blue line) and (e)-(f) (red line). The purple line shows the best reduction that achieved by optimizing $r_i$ and $a_i$ fault by fault.
The inset highlights the small discrepancies induced by the choice of $\eta_b$ for the faults
with largest impact.}\label{fig:vulnerability}
\end{figure}

At the end of the previous section
we argued that an homogeneous distribution of inertia, together with primary control 
increased on the slowest eigenmodes of the network Laplacian minimize
the global vulnerability measure $\mathcal{V}$ of \eqref{eq:vulnerability}
for $\eta_b\equiv1$. This conclusion is confirmed numerically in Fig.~\ref{fig:vulnerability} (a) and (b). 
The optimal placement of primary control displayed in panel (b) decreases $\mathcal{V}$ by more than 12\% with respect to the homogeneous case. 

Setting $\eta_b\equiv 1$ in \eqref{eq:vulnerability} is convenient mathematically, however it treats
all faults equally, regardless of their impact. One may instead adapt $\eta_b$ to obtain inertia
and primary control distributions that reduce the impact of the strongest faults with largest 
$\mathcal{M}_b$. We do this in two different ways, first with
$\eta_b=\mathcal{M}_b^{(0)2}$ and second with 
\begin{equation}
\eta_b = \left\{\begin{array}{l}1\,,\text{ if $\mathcal{M}_b^{(0)}>\mathcal{M}_{\rm thres}$,}\\
0\,,\text{ otherwise.}\end{array}\right. \label{eq:thres}
\end{equation}
The corresponding geographical distributions of inertia and primary control redistribution parameters
$r_i$ and $a_i$ parameters are shown in Fig.~\ref{fig:vulnerability} (c) and (d) and 
Fig.~\ref{fig:vulnerability} (e) and (f) respectively. Compared to the choice $\eta_b\equiv1$ [Fig.~\ref{fig:vulnerability} (a) and (b)], we see 
rather small differences. More importantly, the impact of various choices of $\eta_b$ on 
$\mathcal{M}_b$ is almost negligible, as can be seen in Fig.~\ref{fig:vulnerability} (g). In all cases,
our optimization algorithm reduces first and foremost the impact of the strongest faults, with little
or no influence on the faults that have little impact on grid stability.

It is finally interesting to note that our three choices of $\eta_b$ are close to being optimal, especially
when considering the strongest faults. This can be seen in Fig.~\ref{fig:vulnerability} (g), where the 
purple line shows the maximal obtainable reduction, when inertia and primary control distributions 
are optimized individually fault by fault, i.e. with a different redistribution for each fault. The inset
of Fig.~\ref{fig:vulnerability} (g) shows in particular that for the strongest fault, the three considered choices of $\eta_b$ lead to reductions in $\mathcal{M}_b$ that are very close to the maximal
one. We conclude that rather generically, inertia is optimally distributed homogeneously, while 
primary control should be preferentially located on the slow modes of the grid Laplacian.

\section{Conclusion}\label{section:conclusion}

To find the optimal placement of inertia and primary control in electric power grids with limited
such resources is a problem of paramount importance. Here, we have made what we think is an 
important step forward in constructing a perturbative analytical approach to this problem. In this
approach, both inertia and primary control are limited resources, as should be. Most importantly, 
our method goes beyond the usually made assumption of constant inertia to damping ratio. In our approach inertia and primary control can vary spatially independently from one another.
Our results suggest that the optimal inertia distribution is close to 
homogeneous over the whole grid, but that primary control should be reinforced on buses
located on the support of the slower modes of the network Laplacian. Further work should 
try to extend the approach to the next order in perturbation theory. Work along those lines is in
progress. 

\section*{Acknowledgment}
This work was supported by the Swiss National Science Foundation under grants PYAPP2\_154275
and 200020\_182050.

\appendices

\section{}\label{sect:continuation}

\begin{IEEEproof}[Proof of Proposition~\ref{prop:perturbed_xii}]

The proof follows the same steps as for Proposition \ref{proposition1}. The calculations are rather tedious, though algebraically straightforward. In the following, we sketch the 
calculational steps. Assuming that $\bm H$ can be diagonalized as $\bm{t}^{R}\,\bm \mu \,\bm{t}^{L}$, where  $\bm \mu \equiv {\rm diag}(\{\mu_{\alpha s}\})$ and $\bm{t}^{R,L}$ are matrices containing the right and left eigenvectors of $\bm H$, the problem is resolved by:\\
1) changing variables $\bm \chi\equiv \bm t^{L}[\bm \xi^\top \dot{\bm \xi}^\top]^\top$ to diagonalize \eqref{eq:H1}, as
\begin{equation}
\bm{\dot\chi}=
\bm{\mu}\bm\chi+\bm{t}^L
\left[\!\!\begin{array}{c}
\mathbb{0}_{N\times 1}\\
\bm{\mathcal{P}}
\end{array}\!\!\right]\equiv\bm\mu\bm\chi+\bm{\tilde \mathcal{P}}\,;\label{eq:unperturbed_ev}
\end{equation}
2) solving \eqref{eq:unperturbed_ev} as
\begin{equation}
\chi_{\alpha\pm}=-\frac{\tilde{\mathcal{P}}_{\alpha\pm}}{\mu_{\alpha\pm}}\Big(1-e^{\mu_{\alpha\pm}t}\Big)\,,\;\forall\alpha>1\,;\label{eq:chi_pm}
\end{equation}
3) Obtaining $\dot\xi_\alpha$ via the inverse transformation  $[\bm\xi^\top \bm{\dot\xi}^\top]^\top=\bm t^R\bm\chi$.\\
These three steps are carried out with the approximate expressions $\bm t^{R,L}=\bm t^ {R,L(0)}+g\bm t^{R,L(1)}$ and $\mu_{\alpha\pm}=\mu_{\alpha\pm}^{(0)}+g\mu_{\alpha\pm}^{(1)}$ obtained with the first order in $g$ corrections presented in \eqref{eq:pertval}--\eqref{eq:pertvec}. One gets
\begin{align}
&\!\!\left[\!\!\!\begin{array}{c}\xi_\alpha\\\dot\xi_\alpha\end{array}\!\!\!\right]=\!
\left[\!\begin{array}{cc}1 & 1\\\mu_{\alpha+}^{(0)} & \!\!\!\!\mu_{\alpha-}^{(0)}\end{array}\!\!\right]\!\!\left[\!\!\!\begin{array}{c} \chi_{\alpha+} \\ \chi_{\alpha-}\end{array}\!\!\!\right]
-\frac{g\gamma\bm{V}_{2;\alpha\alpha}}{f_\alpha ^2}\left[\!\begin{array}{cc}\mu_{\alpha+}^{(0)}&\!\!\!\mu_{\alpha-}^{(0)}\\\lambda_\alpha &\!\!\!\lambda_\alpha \end{array}\!\right]\!
\!\!\left[\!\!\!\begin{array}{c} \chi_{\alpha+}^{(0)} \\ \chi_{\alpha-}^{(0)}\end{array}\!\!\!\right]\nonumber\\
&-g\gamma\sum_{\beta\neq \alpha}\frac{\bm{V}_{2;\alpha\beta}}{\lambda_\alpha -\lambda_\beta}
\left[\!\begin{array}{cc} \mu_{\beta+}^{(0)} & \mu_{\beta-}^{(0)} \\ \mu_{\beta+}^{(0)2}   & \mu_{\beta-}^{(0)2}  \end{array}\!\right]\!\!
\left[\!\!\!\begin{array}{c} \chi_{\beta+}^{(0)} \\ \chi_{\beta-}^{(0)}\end{array}\!\!\!\right]+\mathcal{O}(g^2)\,,\label{eq:xi_exp}
\end{align}
with
\begin{align}
\chi_{\alpha\pm}&=-\frac{1}{\mu_{\alpha\pm}^{(0)}}\bigg[\tilde{\mathcal{P}}_{\alpha\pm}^{(0)}+g\tilde{\mathcal{P}}_{\alpha\pm}^{(1)}-g\frac{\mu_{\alpha\pm}^{(1)}\tilde{\mathcal{P}}_{\alpha\pm}^{(0)}}{\mu_{\alpha\pm}^{(0)} }\bigg]\Big(1-e^{\mu_{\alpha\pm}^{(0)}t}\Big)\nonumber\\
&+gt\frac{\mu_{\alpha\pm}^{(1)}\tilde{\mathcal{P}}_{\alpha\pm}^{(0)}}{\mu_{\alpha\pm}^{(0)}} e^{\mu_{\alpha\pm}^{(0)}t}+\mathcal{O}(g^2)\,,\label{eq:chi_exp}
\end{align}
where
\begin{align}
\left[\!\!\!\begin{array}{c}\tilde{\mathcal{P}}_{\alpha+}^{(0)}\\ \tilde{\mathcal{P}}_{\alpha-}^{(0)}\end{array}\!\!\!\right]&=\frac{i}{f_\alpha }
\!\left[\!\!\!\begin{array}{cc}\mu_{\alpha-}^{(0)}& \!\!\!\!-1\\-\mu_{\alpha+}^{(0)}\!\! & 1\end{array}\!\!\!\right]\!
\left[\!\!\begin{array}{c} 0 \\ \mathcal{P}_\alpha \end{array}\!\!\right]\,,\nonumber\\
\left[\!\!\!\begin{array}{c}\tilde{\mathcal{P}}_{\alpha+}^{(1)}\\ \tilde{\mathcal{P}}_{\alpha-}^{(1)}\end{array}\!\!\!\right]
&=\frac{i\gamma}{f_\alpha}
\bigg(-\frac{\bm{V}_{2;\alpha\alpha}}{f_\alpha^2}\left[\!\begin{array}{cc}\lambda_\alpha &\!\!\!\!\!\!-\mu_{\alpha-}^{(0)}\\\!\!\!-\lambda_\alpha &\mu_{\alpha+}^{(0)}\end{array}\!\!\!\right]\!\!
\left[\!\!\begin{array}{c} 0 \\ \mathcal{P}_\alpha \end{array}\!\!\right]\nonumber\\
&+\sum_{\beta\neq \alpha}\frac{\bm{V}_{2;\alpha\beta}}{(\lambda_\alpha -\lambda_\beta)}\!\!
\left[\!\begin{array}{cc}\lambda_\beta &\!\!\!-\mu_{\alpha+}^{(0)}\\\!-\lambda_\beta &\mu_{\alpha-}^{(0)}\end{array}\!\!\!\right]\!\!\left[\!\begin{array}{c} 0 \\  \mathcal{P}_\beta\end{array}\!\right]\bigg)\,.\label{eq:P}
\end{align}
\eqref{eq:perturbed_xii} is obtained from \eqref{eq:xi_exp} by applying trigonometric identities.
\end{IEEEproof}

\begin{IEEEproof}[Proof of Theorem~\ref{thm:lin_opt_inertia}] To leading order in $\mu=\delta m/m$, this optimization problem is equivalent to the following linear programming problem \cite{bertsimas1997introduction} 
\begin{align}
&\min_{\{r_i\}}\sum_i\rho_ir_i\,,\\
\text{s.t.}&\;|r_i|\le 1\,,\\
&\sum_i r_i=0\,.
\end{align}
It is solved by the Lagrange multipliers method, with the Lagrangian function
\begin{equation}
\mathcal{L}=\sum_{i=1}^N\rho_ir_i+\sum_{i=1}^N\varepsilon_i(r_i^2-1)+\varepsilon_0\sum_{i=1}^Nr_i\,,
\end{equation}
where $\varepsilon_i$ and $\varepsilon_0$ are Lagrange multipliers. We get
\begin{equation}
\frac{\partial\mathcal{L}}{\partial r_i}=\rho_i+2\varepsilon_ir_i+\varepsilon_0=0\,,\;\forall i\,.\label{eq:lagrange}
\end{equation}
The solution must satisfy the Karush-Kuhn-Tucker (KKT) conditions \cite{bertsimas1997introduction}, in particular the complementary slackness (CS) condition which imposes that either $\varepsilon_i=0$ or $r_i=\pm1\,,\;\forall i$. The former choice leads generally to a contradiction. From \eqref{eq:lagrange} and dual feasibility condition, one gets
\begin{equation}
\varepsilon_i=-\frac{\varepsilon_0+\rho_i}{2r_i}\ge 0\,.\label{eq:dual}
\end{equation}
This imposes that $r_i=-{\rm sgn}(\varepsilon_0+\rho_i)$. To ensure that $\sum_ir_i=0$ is satisfied, $\varepsilon_0$  is set to minus the median value of $\rho_i$. If the number of bus $N$ is odd, the $r_i$ corresponding to the median value of $\rho_i$ is set to zero.
\end{IEEEproof}

\begin{thebibliography}{10}
\providecommand{\url}[1]{#1}
\csname url@samestyle\endcsname
\providecommand{\newblock}{\relax}
\providecommand{\bibinfo}[2]{#2}
\providecommand{\BIBentrySTDinterwordspacing}{\spaceskip=0pt\relax}
\providecommand{\BIBentryALTinterwordstretchfactor}{4}
\providecommand{\BIBentryALTinterwordspacing}{\spaceskip=\fontdimen2\font plus
\BIBentryALTinterwordstretchfactor\fontdimen3\font minus
  \fontdimen4\font\relax}
\providecommand{\BIBforeignlanguage}[2]{{%
\expandafter\ifx\csname l@#1\endcsname\relax
\typeout{** WARNING: IEEEtran.bst: No hyphenation pattern has been}%
\typeout{** loaded for the language `#1'. Using the pattern for}%
\typeout{** the default language instead.}%
\else
\language=\csname l@#1\endcsname
\fi
#2}}
\providecommand{\BIBdecl}{\relax}
\BIBdecl

\bibitem{ulbig2014impact}
A.~Ulbig, T.~S. Borsche, and G.~Andersson, ``Impact of low rotational inertia
  on power system stability and operation,'' \emph{IFAC Proceedings Volumes},
  vol.~47, no.~3, pp. 7290--7297, 2014.

\bibitem{Mil15}
M.~Milligan, B.~Frew, B.~Kirby, M.~Schuerger, K.~Clark, D.~Lew, P.~Denholm,
  B.~Zavadil, M.~O'Malley, and B.~Tsuchida, ``Alternatives no more: Wind and
  so- lar power are mainstays of a clean, reliable, affordable grid,''
  \emph{IEEE Power and Energy Magazine}, no.~13, pp. 78--87, 2015.

\bibitem{Win15}
W.~Winter, K.~Elkington, G.~Bareux, and J.~Kostevc, ``Pushing the limits:
  Europe's new grid: Innovative tools to combat transmission bottlenecks and
  reduced inertia,'' \emph{IEEE Power and Energy Magazine}, no.~13, pp. 60--74,
  2015.

\bibitem{borsche2015effects}
T.~S. Borsche, T.~Liu, and D.~J. Hill, ``Effects of rotational inertia on power
  system damping and frequency transients,'' in \emph{Decision and Control
  (CDC), 2015 IEEE 54th Annual Conference on}.\hskip 1em plus 0.5em minus
  0.4em\relax IEEE, 2015, pp. 5940--5946.

\bibitem{mevsanovic2016comparison}
A.~Me{\v{s}}anovi{\'c}, U.~M{\"u}nz, and C.~Heyde, ``{Comparison of
  $H_{\infty}$, $H_2$, and pole optimization for power system oscillation
  damping with remote renewable generation},'' \emph{IFAC-PapersOnLine},
  vol.~49, no.~27, pp. 103--108, 2016.

\bibitem{poolla2017optimal}
B.~K. Poolla, S.~Bolognani, and F.~D{\"o}rfler, ``Optimal placement of virtual
  inertia in power grids,'' \emph{IEEE Transactions on Automatic Control},
  vol.~62, no.~12, pp. 6209--6220, 2017.

\bibitem{paganini2017global}
F.~Paganini and E.~Mallada, ``Global performance metrics for synchronization of
  heterogeneously rated power systems: The role of machine models and
  inertia,'' in \emph{55th Annual Allerton Conference on Communication,
  Control, and Computing}, Oct 2017, pp. 324--331.

\bibitem{bor18}
T.~S. Borsche and D{\"o}rfler, ``On placement of synthetic inertia with
  explicit time-domain constraints,'' \emph{arXiv:1705.03244}, 2017.

\bibitem{pagnier2019inertia}
L.~Pagnier and P.~Jacquod, ``Inertia location and slow network modes determine
  disturbance propagation in large-scale power grids,'' \emph{PLoS ONE},
  vol.~14, no.~3, p. e0213550, 2019.

\bibitem{stewart1990matrix}
G.~W. Stewart and J.-G. Sun, \emph{Matrix perturbation theory}.\hskip 1em plus
  0.5em minus 0.4em\relax Academic Press, Boston, 1990.

\bibitem{macintyre}
D.~McIntyre, \emph{Quantum Mechanics}.\hskip 1em plus 0.5em minus 0.4em\relax
  Pearson Addison-Wesley, San Francisco, 2012.

\bibitem{coletta2018performance}
T.~Coletta and P.~Jacquod, ``Performance measures in electric power networks
  under line contingencies,'' \emph{IFAC PapersOnLine}, no. 51--23, pp.
  337--342, 2018.

\bibitem{Pir17}
M.~Pirani, J.~Simpson-Porco, and B.~Fidan, ``System-theoretic performance
  metrics for low-inertia stability of power networks,'' in \emph{Decision and
  Control (CDC), 2017 IEEE 56th Annual Conference on}.\hskip 1em plus 0.5em
  minus 0.4em\relax IEEE, 2017, pp. 5106--5111.

\bibitem{Guo18}
L.~Guo, C.~Zhao, and S.~Low, ``Graph laplacian spectrum and primary frequency
  regulation,'' \emph{arXiv preprint arXiv:1803.03905}, 2018.

\bibitem{Por18}
M.~Porfiri and M.~Frasca, ``Robustness of synchronization to additivie noise:
  how vulnerability depends on dynamics,'' \emph{to appear in IEEE}, 2018.

\bibitem{machowski2008power}
J.~Machowski, J.~Bialek, and J.~R. Bumby, \emph{Power system dynamics:
  stability and control}, 2nd~ed.\hskip 1em plus 0.5em minus 0.4em\relax John
  Wiley \& Sons, 2008.

\bibitem{coletta2018transienta}
T.~Coletta, B.~Bamieh, and P.~Jacquod, ``Transient performance of electric
  power networks under colored noise,'' in \emph{2018 IEEE Conference on
  Decision and Control (CDC)}.\hskip 1em plus 0.5em minus 0.4em\relax IEEE,
  2018, pp. 6163--6167.

\bibitem{tegling2015price}
E.~Tegling, B.~Bamieh, and D.~F. Gayme, ``The price of synchrony: Evaluating
  the resistive losses in synchronizing power networks.'' \emph{IEEE Trans.
  Control of Network Systems}, vol.~2, no.~3, pp. 254--266, 2015.

\bibitem{Fardad14}
M.~Fardad, F.~Lin, and M.~R. Jovanovic, ``Design of optimal sparse
  interconnection graphs for synchronization of oscillator networks,''
  \emph{IEEE Transactions on Automatic Control}, vol.~59, no.~9, pp.
  2457--2462, 2014.

\bibitem{gayme16}
T.~W. Grunberg and D.~F. Gayme, ``Performance measures for linear oscillator
  networks over arbitrary graphs,'' \emph{IEEE Transactions on Control of
  Network Systems}, vol.~PP, no.~99, p.~1, 2016.

\bibitem{siami2016fundamental}
M.~Siami and N.~Motee, ``Fundamental limits and tradeoffs on disturbance
  propagation in linear dynamical networks,'' \emph{IEEE Transactions on
  Automatic Control}, vol.~61, no.~12, pp. 4055--4062, 2016.

\bibitem{tyloo2018robustness}
M.~Tyloo, T.~Coletta, and P.~Jacquod, ``Robustness of synchrony in complex
  networks and generalized kirchhoff indices,'' \emph{Phys. Rev. Lett.}, vol.
  120, p. 084101, 2018.

\bibitem{guo2019performance}
Y.~Guo and T.~H. Summers, ``A performance and stability analysis of low-inertia
  power grids with stochastic system inertia,'' \emph{arXiv preprint
  arXiv:1903.00635}, 2019.

\bibitem{Kle93}
D.~Klein and M.~Randi\'c, ``Resistance distance,'' \emph{J. Math. Chem.},
  vol.~12, p.~81, 1993.

\bibitem{tyloo2018key}
M.~Tyloo, L.~Pagnier, and P.~Jacquod, ``The key player problem in complex
  oscillator networks and electric power grids: Resistance centralities
  identify local vulnerabilities,'' submitted.

\bibitem{Gut96}
I.~Gutman and B.~Mohar, ``The quasi-wiener and the kirchhoff indices
  coincide,'' \emph{J. Chem. Inf. Comput. Sci.}, vol.~36, no.~2, p. 982, 1996.

\bibitem{Zhu96}
H.~Zhu, D.~Klein, and I.~Lukovits, ``Extensions of the wiener number,''
  \emph{J. Chem. Inf. Comput. Sci.}, vol.~36, no.~2, p. 420, 1996.

\bibitem{Gam17}
L.~Gambuzza, A.~Buscarino, L.~Fortuna, M.~Porfiri, and M.~Frasca, ``Analysis of
  dynamical robustness to noise in power grids,'' \emph{IEEE Journal on
  Emerging and Selected topics in Circuits and Systems}, no.~7, p. 413, 2017.

\bibitem{dorfler2013kron}
F.~D{\"o}rfler and F.~Bullo, ``Kron reduction of graphs with applications to
  electrical networks.'' \emph{IEEE Trans. on Circuits and Systems}, vol.~60,
  no.~1, pp. 150--163, 2013.

\bibitem{bertsimas1997introduction}
D.~Bertsimas and J.~N. Tsitsiklis, \emph{Introduction to linear
  optimization}.\hskip 1em plus 0.5em minus 0.4em\relax Athena Scientific
  Belmont, MA, 1997.

\end{thebibliography}

\begin{IEEEbiography}[{\includegraphics[width=1in,height=1.25in,clip,keepaspectratio]{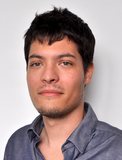}}]{Laurent Pagnier} (M'17) received the M.S. and Ph.D. degrees in theoretical physics from the EPFL, Lausanne, Switzerland, in 2014 and 2019 respectively. He is currently working as a postdoctoral researcher with the Electrical Energy Efficiency Group, University of Applied Sciences of Western Switzerland, Sion, Switzerland. His research interests include the development of new renewable energy sources and its effects on transmission grids.
\end{IEEEbiography}

 \begin{IEEEbiography}[{\includegraphics[width=1in,height=1.25in,clip,keepaspectratio]{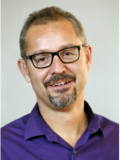}}]{Philippe Jacquod} (M'16) received the Diplom degree in theoretical physics from the ETHZ, Z\"{u}rich, Switzerland, in 1992, and the PhD degree in
 natural sciences from the University of Neuch\^{a}tel,
 Switzerland, in 1997. He is a professor with the engineering department,
 University of Applied Sciences of Western Switzerland, Sion, Switzerland, with a joint appointment with the Department of Quantum Matter Physics, University of Geneva, Switzerland. From 2003 to 2005 he was an assistant professor with the
 theoretical physics department, University of Geneva,
 Switzerland and from 2005 to 2013 he was a professor
 with the physics department, University of Arizona, Tucson, USA. His main research topics is in power systems and how they
 evolve as the energy transition unfolds. He has published about 100
 papers in international journals, books and conference proceedings.
 \end{IEEEbiography}

\end{document}